\documentclass[]{amsart}
\usepackage{amsmath,amssymb,amsthm,mathptmx,dsfont,framed, graphicx, textcomp,mathtools}
\usepackage{tikz, braids}
\usepackage[T1]{fontenc} 
\usepackage{enumerate}    
\usepackage[hidelinks]{hyperref}
\usepackage{appendix}
\usepackage{epsfig}
\usepackage{verbatim}
\usepackage{bm}
\usepackage{mathrsfs}
\usepackage{yfonts}
\usepackage{float}
\usepackage{times}
\usepackage{amsfonts}
\usepackage{pst-all}
\usepackage{leftidx}

\usetikzlibrary{shapes,snakes,positioning,decorations.pathmorphing,calc,backgrounds,chains,fit,matrix,trees,scopes,shapes.geometric,shapes.multipart}
\usetikzlibrary{decorations.pathreplacing, arrows}
\usetikzlibrary{decorations.markings}

\newcommand{\mymk}[1]{%
  \tikz[baseline=(char.base)]\node[thick, black, anchor=south west, draw,ellipse, inner sep=2pt, minimum size=7mm,
    text height=2mm](char){\textcolor{black}{#1}} ;}

\pgfdeclarelayer{top}
\pgfsetlayers{main,top}

\newtheorem{thm}{Theorem}[section]

\newtheorem{cor}[thm]{Corollary}
\newtheorem{conj}[thm]{Conjecture}
\newtheorem{defn}[thm]{Definition}

\DeclareMathAlphabet{\mathcal}{OMS}{cmsy}{m}{n}
\DeclareMathOperator{\Hom}{Hom} 
 
\DeclareMathOperator{\Span}{Span} 

\newcommand{\ZZ}{\mathds{Z}}
\DeclareMathOperator{\Aut}{Aut}
\DeclareMathOperator{\Rep}{Rep}
\DeclareMathOperator{\Obj}{Obj}
\usepackage[all]{xy}

\newcommand{\BB}{\mathcal{B}}
\newcommand{\CC}{\mathcal{C}}
\newcommand{\DD}{\mathcal{D}}

\newcommand{\btikz}{\begin{tikzpicture}}
\newcommand{\etikz}{\end{tikzpicture}}
\newcommand{\AutC}{\text{Aut}_{br}^{\otimes}(\CC)}
\newcommand{\gammabar}{\ensuremath\gamma\kern-0.53em-}

\newcommand{\stdlabels}[4]{\draw (1.5,0) node[below] {$\phantom{X}\sigma\phantom{X}$};

\draw (2.5,0) node[below] {$\phantom{X}#1\phantom{X}$};
\draw (3.5,0) node[below] {$\phantom{X}\sigma\phantom{X}$};
\draw (4.5,0) node[below] {$1$};
\draw (5.5,0) node[below] {$X$};
\draw (6.5,0) node[below] {$\phantom{X}#2\phantom{X}$};
\draw (7.5,0) node[below] {$X$};
\draw (1.5,10) node[above] {$\phantom{X}\sigma\phantom{X}$};
\draw (2.5,10) node[above] {$\phantom{X}#3\phantom{X}$};
\draw (3.5,10) node[above] {$\phantom{X}\sigma\phantom{X}$};
\draw (4.5,10) node[above] {$1$};
\draw (5.5,10) node[above] {$X$};
\draw (6.5,10) node[above] {$\phantom{X}#4\phantom{X}$};
\draw (7.5,10) node[above] {$X$};}

\newcommand{\labelstwo}[2]{\draw (1.5,0) node[below] {$\phantom{X}\sigma\phantom{X}$};
\draw (3.5,0) node[below] {$\phantom{X}#1\phantom{X}$};
\draw (5.5,0) node[below] {$X$};
\draw (6.5,0) node[below] {$\phantom{X}1\phantom{X}$};
\draw (7.5,0) node[below] {$X$};
\draw (1.5,10) node[above] {$\phantom{X}\sigma\phantom{X}$};
\draw (3.5,10) node[above] {$\phantom{X}#2\phantom{X}$};
\draw (5.5,10) node[above] {$X$};
\draw (6.5,10) node[above] {$\phantom{X}1\phantom{X}$};
\draw (7.5,10) node[above] {$X$};}

\newcommand{\labelsthree}{\draw (1.5,0) node[below] {$\phantom{X}\sigma\phantom{X}$};
\draw (3.5,0) node[below] {$\phantom{X}1\phantom{X}$};
\draw (6,0) node[below] {$X$};
\draw (7,0) node[below] {$X$};
\draw (1.5,10) node[above] {$\phantom{X}\sigma\phantom{X}$};
\draw (3.5,10) node[above] {$\phantom{X}1\phantom{X}$};
\draw (6,10) node[above] {$X$};
\draw (7,10) node[above] {$X$};}

\begin{document}

\title{Symmetry defects and their application to topological quantum computing}

\author{Colleen Delaney}
\email{cdelaney@math.ucsb.edu}
\address{Department of Mathematics,
    University of California,
    Santa Barbara, CA
    U.S.A.}

\author{Zhenghan Wang}
\email{zhenghwa@microsoft.com}
\address{Microsoft Research Station Q and Department of Mathematics,
    University of California,
    Santa Barbara, CA
    U.S.A.}
\maketitle

\begin{abstract}
We describe the mathematical theory of topological quantum computing with symmetry defects in the language of fusion categories and unitary representations. Symmetry defects together with anyons are modeled by $G$-crossed braided extensions of unitary modular tensor categories. The algebraic data of these categories afford projective unitary representations of the braid group.  Elements in the image of such representations correspond to quantum gates arising from exchanging anyons and symmetry defects in topological phases of matter with symmetry. We provide some small examples that highlight features of practical interest for quantum computing. In particular, symmetry defects show the potential to generate non-abelian statistics from abelian topological phases and to be used as as a tool to enlarge the set of quantum gates accessible to an anyonic device. 
\end{abstract}

\section{Introduction}
Topological order of a 2-dimensional topological phase of matter can be modeled by a unitary modular tensor category $\mathcal{C}$.  If a 2-dimensional topological phase of matter has a symmetry group $G$, then $\mathcal{C}$ inherits an action of $G$.  A nontrivial action of $G$ on $\mathcal{C}$ is a form of symmetry breaking, called weak symmetry breaking in \cite{K1}. It follows that the associated topological defects, which manifest the topological properties of symmetry breaking patterns, provide a window to the world of symmetry and topology.  A practical reason to study such symmetry defects is as a new mechanism to generate non-abelian statistics from a purely abelian topological phase.  For example, we show that the symmetry defects arising from the electromagnetic duality of the $\mathds{Z}_2$ toric code behave almost identically to Majorana zero modes or Ising anyons.

The symmetry defects we consider here are extrinsic point-like objects that are in general confined.  Like anyons, they lead to topological degeneracy and braid group representations, though only projectively.  In this paper, we investigate the mathematical modeling of symmetry defects and their application to topological quantum computing.

In Section 2 we review the algebraic description of symmetry defects via $G$-crossed extensions of unitary modular tensor categories following \cite{BBCW}, where this mathematical model for topological order with symmetry was first introduced. Some examples are presented to illustrate the theory and are further developed in Section 3, where we describe the (projective) unitary representations associated to unitary $G$-crossed extensions and discuss their interpretation as quantum gates for a topological quantum computer built by braiding symmetry defects and anyons. In Section 4, we give a mathematical analysis of an example of topological order with $\mathds{Z}_2$ symmetry that shows  allowing projective measurement and braiding with symmetry defects can enlarge an anyonic gate set with finite image in the unitary group to one that is universal for quantum computation. This feature of two layers of the Ising theory with a layer-exchange symmetry $\left ( \textbf{Ising}^{(1)} \boxtimes \textbf{Ising}^{(1)}\right)^{\times}_{\mathds{Z}_2}$ was proposed in \cite{BJQ} based on hypothetical data for the $\mathds{Z}_2$-crossed extension and a physical argument. We conclude in Section 5 with some comments about future directions, physical and mathematical.

\section{Algebraic theory of symmetry defects}
Symmetry defects are point-like objects that are not intrinsic quasiparticle excitations of a topological phase of matter, but can nevertheless be used for topological quantum computing (TQC). In analogy with the mathematical definition of an anyon as a simple object in a unitary modular tensor category (UMTC), a symmetry defect can be defined as a simple object in a non-trivial sector of a unitary $G$-crossed extension $\CC^{\times}_G$ of a UMTC $\CC$ \cite{BBCW}. The mathematical model of symmetry defects and anyons is then a \emph{unitary $G$-crossed braided fusion category}, UGxBFC for short. Since in physics symmetry defect models arise from anyon models $\CC$ with symmetry $G$, we begin with the construction of UGxBFCs as extensions of UMTCs and delay an abstract definition of a UGxBFC until later in the section. Both notions will be useful to have on hand when we discuss how symmetry defect models can also come from a topological phase $\mathcal{D}$ in a different way: as the first step in a two-step ``topological Bose condensation'' procedure (de-equivariantization in the language of fusion categories) \cite{EGNO}. 

\subsection{$G$-crossed extensions of UMTCs}
Here we briefly recall the construction of a $G$-crossed category from a UMTC $\CC$ and identify symmetry defects and their mathematical properties. Our discussion follows \cite{BBCW} but our notation is slightly different in that we use standard font for group elements.

\subsubsection{Preliminaries}

In what follows we assume $G$ is a finite group. Although $G$ can be non-abelian we will consistently use $0$ to denote the identity element. We will work with a certain ``categorification'' of $G$, a categorical group $\mathcal{G}$.  The objects of $\mathcal{G}$ are group elements and morphism sets $\Hom(g,h)=\delta_{gh}$ of two objects $g$ and $h$ are empty if $g\ne h$ and contain only the identity if $g=h$. That is, there are no nontrivial morphisms between distinct group elements, and the only morphism between identical group elements is the identity. With respect to the tensor product given by group multiplication $g \otimes h = gh$, $\mathcal{G}$ forms a strict monoidal category. 

We will use $\CC$ to denote an arbitrary UMTC. We denote the set of objects of $\CC$ by $\text{Obj}(\CC)$. Recall that an anyon is a simple object in $\CC$. We call the isomorphism class of an anyon its type. When we work with concrete examples we always pass to a skeletonization of $\CC$, in which case there is only one object per isomorphism class and we conflate anyons and anyon types. We call an anyon abelian if it has quantum dimension one. This definition is equivalent to the others in use, for example that in terms of images of associated braid group representations by the work of \cite{RW}.  The subset of $\Obj(\CC)$ consisting of abelian anyon types forms a subgroup under the tensor product, which we denote by $\mathcal{A}$. To see why this is true, observe that the tensor unit (the vacuum) has quantum dimension one and the quantum dimension function $d_{-}: \Obj(\CC) \to \mathbb{R}_{\ge 0}$ is linear and multiplicative, in the sense that
$$a \otimes b = \bigoplus_c N^{ab}_c c \implies d_a\cdot d_b = \sum_{c} N^{ab}_c d_c.$$
In particular, if $d_a=d_b=1$, then $a \otimes b = c$ where $d_c=1$. It also follows that the dual $a^*$ of an abelian anyon $a$ is also abelian, since $a \otimes a^* = 1$ implies $d_a^*=1$.  It follows that the subset of abelian anyon types contains the vacuum and is closed under multiplication and inverses. Since fusion is associative on the level of isomorphism classes, $\mathcal{A}$ forms a group. 

We denote by $\AutC$ the group of braided-tensor autoequivalences of $\CC$ up to natural isomorphism under composition. We will also need the categorical group $\underline{\AutC}$, with group elements of $\AutC$, tensor product given by composition, and morphisms given by natural isomorphisms of functors.  

\subsubsection{From group symmetry to group extension theory of MTCs}
To construct a UGxBFC, we start with a UMTC $\CC$ and a group homomorphism $\rho: G \to \AutC$, called a \emph{global symmetry}. When this group homomorphism $\rho$ can be promoted to a monoidal functor $\underline{\rho}:\mathcal{G} \to \underline{\AutC}$, it is called a \emph{categorical global symmetry} of $\CC$. The existence of such a lift $\underline{\rho}$ depends on the vanishing of a certain cohomological obstruction in $H^3(G, \mathcal{A})$, where $\mathcal{A}$ is the subgroup of abelian anyons under fusion. When this obstruction vanishes, the different extensions are ``classified" by  $H^2_{\rho}(G,\mathcal{A})$ in the sense that the distinct liftings of the global symmetry $\rho$ to the categorical global symmetry $\underline{\rho}$ form a torsor over $H^2_{\rho}(G,\mathcal{A})$.

With a categorical global symmetry in hand, the next ingredient for constructing a $G$-extension of $\CC$ is an isomorphism of the categorical groups $\underline{\AutC}$ and $\underline{\text{Pic}(\CC)}$, the categorical Picard group consisting of invertible bimodule categories of $\CC$ \cite{ENO}. This associates to each categorical braided tensor autoequivalence $\underline{\rho_g}$ an invertible $\CC$-bimodule category $\CC_g$, a correspondence which is at the heart of the relationship between symmetry in the abstract and symmetry defects. When an additional cohomological obstruction in $H^4(G;U(1))$ vanishes, consistent fusion of invertible bimodule categories $\CC_g \boxtimes \CC_h$ can be defined so that elements of $\CC_g$ satisfy the pentagon equations. The different solutions to the pentagon equations then form a torsor over $H^3(G,U(1))$. 
 
Each such solution corresponds noncanonically to a unitary fusion category denoted by $\CC^{\times}_G=\bigoplus_{g \in G} \CC_g$, where $\CC_0=\CC$. The inequivalent $G$-crossed braided extensions of a UMTC $\CC$ are classified by the data $(\rho, \alpha, \beta)$, where $\alpha \in H^2_{\rho}(G,\mathcal{A})$ and $\beta \in H^3(G,U(1))$ \cite{BBCW, ENO}. We will see in the next subsection that $\CC^{\times}_G$ has the structure of a unitary \emph{$G$-crossed braided fusion category}.

Now the anyons of $\CC$ are identified as simple objects in $\CC_0$ and the symmetry defects as simple objects in $\bigoplus_{g\ne 0}\CC_g$:

\begin{defn} A symmetry defect with flux $g$ is a simple object in the invertible $\CC_0$-bimodule category $\CC_g$. 
\end{defn}
We will sometimes call the bimodule category $\CC_g$ the $g$-flux sector, or a $g$-sector.

\subsubsection{Physical interpretation of UGxBFCs}
We note that in the physics literature, especially in the early development of TQC, the word ``anyon" was used in a broad sense to describe both intrinsic topological order and extrinsic defects \cite{DSFN}. Here we use the term ``anyon" for simple objects in a UMTC, and reserve the term ``symmetry defects" for simple objects in a nontrivial sector of a UGxBFC. While in this framework anyons could be considered as trivial defects, the physics of extrinsic defects differs in important ways from quasiparticle excitations of a topological phase of matter and justifies a distinction in the terminology.

The mathematical formalism describing symmetry defects, namely unitary $G$-crossed extensions of UMTCs, holds for any finite group $G$ acting by a categorical global symmetry (when the $H^4$ obstruction vanishes). However, there are further considerations to make about the nature of the categorical group action $\underline{\rho}$ when trying to identify them as models for quantum systems. Apart from unitarity, generally one is interested in symmetries that are local in the following sense: operators that effect the symmetry transformations can be written as a product of local operators that remain localized in the same region under the symmetry action. These are called ``on-site'' symmetries. When the group action is unitary and on-site, the resulting class of UGxBFCs characterize what are known in condensed matter theory as (2+1)D symmetry-enriched topological (SET) phases of matter. However, not all physically meaningful symmetries are unitary and ``on-site'': the authors of \cite{BBCW} make sense of anti-unitary and ``quasi-on-site'' symmetries, which can include time-reversal, translational, rotational, and other spacetime symmetries.  

\subsection{$G$-crossed braided fusion categories}
We will see that a $G$xBFC is a fusion category together with three compatible structures: a $G$-action, a $G$-grading, and a $G$-braiding. Below $\CC$ will denote a fusion category, not a UMTC as in the previous subsection.

Recall from the previous subsection that a $G$-action is a monoidal functor $\rho: \mathcal{G} \to \underline{\AutC}$. We write the $G$-action on objects of $\CC$ as $\rho_g(X)=\prescript{g}{}{X}$. 

In the fusion category setting, a $G$-grading is a map $\partial: \text{Obj}( \CC_{hom}) \to G$, where $\CC_{hom}$ is a full fusion subcategory of $\CC$ with the property that it generates the objects of $\CC$ by direct sum \cite{M}. The full subcategories $\CC_g = \partial^{-1}(g)$ given by the $G$-grading are the $g$-flux sectors seen previously.

\begin{defn}A $G$-crossed braided fusion category $\CC^{\times}_G$ is a fusion category with 

\begin{itemize} 
\item a $G$-action $^g(\cdot): \mathcal{G}\to \Aut_{\otimes}(\CC^{\times}_G)$
\item a $G$-grading $\partial$ \\
such that $\partial(^gX)=g \left(\partial X\right) g^{-1}$ for all $X \in C_{hom}$
\item a $G$-crossed braiding:  \\
 natural isomorphisms $\{c_{X,Y}: X \otimes Y \to  \prescript{g}{}{Y} \otimes X\}$ for all $X \in \CC_g$, $Y \in \CC_0$, $g \in G$
 \end{itemize}
  satisfying compatibility conditions with one another and the associators \cite{EGNO}, \cite{M}.

\end{defn}

The following theorem from \cite{BBCW}, adapted from \cite{ENO}, says that the $G$-crossed extension of a UMTC from Section 2.1 satisfies the above definition. 

\begin{thm} The unitary $G$-crossed extension $\CC^{\times}_G$ of a UMTC $\CC$ has a canonical $G$-braiding and categorical $G$-action that make $\CC^{\times}_G$ into a unitary $G$-crossed braided fusion category. 
\end{thm}

In particular, the $G$-action on a unitary $G$-crossed extension of a UMTC $\CC$ comes from extending the categorical action from $\CC_0$ to all of the sectors $\CC_g$ and the $G$-grading satisfies the property that the action of $g$ on an object in the sector $\CC_h$ sends it to an object in the sector $\CC_{ghg^{-1}}$.  

To describe topological phases of matter with symmetry we needed $G$-crossed braided fusion categories which were unitary. We have omitted the abstract definition of unitarity in fusion categories and instead will simply state what it means for a skeletonization of a UMTC or UGxBGC to be unitary. This is sufficient for our purposes since in order to work out concrete elements in the image of the projective braid group representations associated to a GxBFC $\CC^{\times}_G$ one must pass to a skeletonization of $\CC^{\times}_G$. Next we discuss these skeletonizations, which we refer to as working with the algebraic data of a UMTC or UGxBFC.

\subsubsection{Algebraic data of a $G$xBFC }

It is well established that a UMTC is characterized by the algebraic data defining its skeletonization, and can effectively be defined as a collection of numbers $\{N_{ab}^c, [F^{abc}_d]_{ef}, R_{c}^{ab}\}$ satisfying certain compatibility conditions. Then an MTC is unitary if there is a basis in which the matrices $[F^{abc}_d]$ and $R_c$ are unitary. Unitarity of the $R$- and $F$- symbols then implies that the associated braid group representations are unitary, endowing them with an interpretation as quantum mechanical processes. 

In \cite{BBCW} this result is extended to GxBFCs, and hence it is justified to define a UGxBFC as a collection of numbers $\{N_{ab}^c, [F^{abc}_d]_{ef}, R_{c}^{ab}, U_{k}(a,b;c), \eta_x(g,h)\}$ which form unitary matrices and satisfy various consistency equations. Here the $R$-symbols keep track of the $G$-crossed braiding, which we recall is not a braiding in the usual sense like it is in a MTC. 

When the topological charges are all anyons (labeled by the trivial element of $G$), the corresponding symbols will obey the usual consistency equations for a UMTC, for example the hexagon equation that encodes compatibility of the associator and braiding in $\CC_0$. The symbols carrying symmetry defect labels now obey more complicated $G$-crossed consistency equations. There is an enlarged set of pentagon equations that must be satisfied by the $F$-symbols, and now compatibility of the associator and $G$-crossed braiding is given by a ``heptagon'' equation see Section 2.4. The additional symbols $U$ and $\eta$ arise from the categorical symmetry, and give rise to additional compatibility relations that the symbols must satisfy. Restricted to the $C_0$ sector, the $U$-symbols encode how the symmetry acts on morphism spaces and the $eta$-symbols encode the natural isomorphism identifying $\rho_{gh}$ with $\rho_g\rho_h$ \cite{BBCW}.

\subsubsection{Computational methods for finding algebraic data}
If one has an anyon model with symmetry, the direct way to determine the $G$-crossed categorical data is to look for unitary solutions to the $G$-crossed consistency equations. This is how the algebraic data for the examples in Section 3 were found \cite{BBCW}. Of course in general it is a difficult problem to solve the consistency equations for a UMTC, let alone for a $G$-extension. 

Another approach to deriving the algebraic data for a symmetry defect model is to make use of ``two-step" gauging \cite{CGPW}. For a UMTC $\CC$ with a symmetry $G$, there is a two-part procedure to gauge the symmetry $G$ and produce another UMTC $\DD$. The first step in this process is to take a $G$-crossed extension $\CC^{\times}_G$ of $\CC$; the second is to $G$-equivariantize $\CC^{\times}_G$ to get the category $\DD=\left(\CC^{\times}_G \right)^G$. Equivariantization is a functor of fusion categories for which there is an inverse functor, called de-equivariantization. De-equivariantization is the process of passing from a fusion category $\DD$ containing $\Rep(G)$ as a symmetric fusion subcategory to the category of $A$-modules over $\DD$, where $A=\text{Fun}(G,\mathds{C})$. It is known that equivariantization and de-equivariantization establish a bijection between certain equivalence classes of GxBFCs and MTCs containing $\text{Rep}(G)$ as a symmetric fusion subcategory \cite{EGNO,M}.

\begin{figure}[H]
$$\begin{tikzpicture}[scale=.75, line width=.5*(4/3)]
\draw[->] (-4,.5) node[left] {\Large $\CC$}--(-.25,.5);
\draw (-2,1) node {$G$-crossed extension};
\draw (.25,0) node[above] {\Large $\CC^{\times}_G$};
\draw[->] (0,0) --(0,-4);
\draw (-.5,-3.5) node[left] {equivariantization};
\draw[<-] (.5,0)--(.5,-4);
\draw (5,-.5) node[left] {de-equivariantization};
\draw (-.75,-4) node[below] {\Large$\Rep(G) \subset \left (\CC^{\times}_G\right)^G$};
\draw (-.75,-5) node[below] {\footnotesize (as a symmetric fusion subcategory)};
\end{tikzpicture}$$
\caption{The relationship between the gauging of an MTC $\CC$ and the GxBFC $\left(\CC\right)^{\times}_G$.}
\end{figure}

As a consequence, if one knows the gauged UMTC $\DD$ or can calculate it from the original UMTC $\CC$, then one can de-equivariantize $\DD$ and calculate the algebraic data for a $G$-crossed theory $\CC^{\times}_G$ explicitly \cite{CHW}. Depending on the category this can be an easier approach to finding solutions to the $G$-crossed consistency equations \cite{DThesis}.

\subsection{Small examples for $G=\mathds{Z}_2$}
A $\ZZ_2$-extension has two graded components, the trivial sector and a nontrivial sector, $\CC^{\times}_{\ZZ_2}=\CC_0 \oplus \CC_1$. It is known that the number of defects in the $g$-flux sector $\CC_g$ is given by the number of fixed points of $\CC_0$ under the action of $g$ \cite{Marcel}. As a consequence one can adopt a convention for labeling $g$-defects in $\CC_g$ by topological charges in $\CC_0$ which are fixed under the action of $g$. We will use this convention with the exception of the toric code symmetry defects, for reasons that will become clear. In Section 3 we will return to these examples and demonstrate how they can be applied for quantum information processing.

\subsubsection{Toric code with electromagnetic duality}
As a UMTC the toric code has four anyon types, given by $\{1,e,m,\psi\}$. The fusion rules are given by the group $\mathds{Z}_2 \times \mathds{Z}_2$ with the identification $1=(0,0),\, e=(1,0),\, m=(0,1),\, \psi=(1,1)$ and have a $\mathds{Z}_2$ symmetry coming from exchanging $e$ and $m$. One can check that the action lifts to a categorical action. Then the $\ZZ_2$-grading is 
$$\text{TC}^{\times}_{\ZZ_2} = \{1,e,m,\psi\}\oplus \{\sigma_+, \sigma_-\},$$
with symmetry defects $\sigma_{\pm}$. All anyons have quantum dimension equal to one, i.e.\ are abelian anyons, while the symmetry defects satisfy $d_{\sigma_{\pm}}=\sqrt{2}$ and have fusion rules with degeneracy.

The table below shows the fusion for the $\mathds{Z}_2$-crossed category. 
$$\begin{array}{c||c|c|c||c|c}

\otimes & e & m & \psi & \sigma_+ & \sigma_-\\
\hline
\hline
e & 1 & \psi & m & \sigma_- & \sigma_+ \\
\hline
m & \psi & 1 & e & \sigma_- & \sigma_+ \\
\hline
\psi & m & e & 1 & \sigma_+ & \sigma_- \\
\hline
\hline
\sigma_+ & \sigma_- & \sigma_- & \sigma_+ & 1 \oplus \psi & e \oplus m \\
\hline
\sigma_- & \sigma_+ & \sigma_+ & \sigma_- & e \oplus m & 1 \oplus \psi \\
\end{array}$$
where we have ommitted the vacuum fusion rules. 

The $\ZZ_2$-crossed $R$-, $F$-, $U$- and $\eta$-symbols were found by solving the $G$-crossed consistency equations in \cite{BBCW}. We will see in the next section that braiding of the symmetry defects leads to non-abelian statistics. 
 
\subsubsection{$\mathds{Z}_3$ anyons with charge conjugation duality}
The UMTC $\ZZ_3$ has three anyon types, labeled by $1, \omega$, and $\omega^*$, with fusion rules given by the group multiplication in $\ZZ_3$. The category has a $\ZZ_2$ symmetry given by interchanging $\omega$ and $\omega^*$, which lifts to a categorical $\ZZ_2$-action and hence admits a $\ZZ_2$-crossed extension which we write

$$(\ZZ_3)^{\times}_{\ZZ_2} = \{1, \omega, \omega^*\} \oplus \{X_{1}\}.$$

As a fusion category, $(\ZZ_3)^{\times}_{\mathds{Z}_2}$ is a Tambara-Yamagami (TY) category, whose algebraic data is known and relatively easy to describe \cite{TY}. Briefly, given a finite group $A$ a TY category has simple objects $A \sqcup \{m\}$ and fusion rules
$$a \otimes b= ab, \hspace{5pt} a \otimes m = m, \hspace{3pt} \text{ and } \hspace{3pt} m \otimes m = \bigoplus_{a \in A} a.$$ A choice of $F$-symbols for TY categories with this fusion are determined by a choice of nondegenerate symmetric bicharacter $\chi:A \times A \to \mathds{C}^{\times}$ and a choice of square root of $|A|$. One can check that in order to correspond to a unitary $\ZZ_2$-crossed BFC the $F$-symbols are determined by a choice of a primitive cube root of unity. For calculations in the next section we choose $\chi$ to take the value $\xi=e^{2\pi i/3}$ on the non-identity diagonal elements of $\ZZ_3 \times \ZZ_3$ and $\bar{\xi}$ on the off-diagonal elements.

For $A=\ZZ_3$ the fusion is given by
$$\begin{array}{c|c|c|c|c}
\otimes & 1 & \omega & \omega^* & X_1 \\
\hline
\hline
1 & 1 & \omega & \omega^* & X_1 \\
\hline
\omega & \omega & \omega^* & 1 & X_1 \\
\hline
\omega^* & \omega^* & 1 & \omega & X_1\\
\hline
X_1 & X_1 & X_1 & X_1 & 1 \oplus \omega \oplus \omega^* \\
\end{array}$$

Formulas for $\mathds{Z}_2$-crossed $R$-symbols were found for general $\ZZ_n$ anyons in \cite{BBCW}.

\subsection{The $G$-crossed graphical calculus}
Before we can apply UGxBFCs for computing we need to introduce the graphical calculus associated to their skeletonizations. The graphical calculus for defects generalizes the graphical calculus for UMTCs in the sense that whenever the topological charges involved in a diagram are all anyons, we recover the graphical calculus of a UMTC. We refer the reader to \cite{BBCW,CBMS} for a review of the diagrammatic calculus for anyons.
As with anyons, the consistency equations that ensure the various categorical structures are compatible with one another are encoded in the topological properties of the diagrams. Here we point out only the parts of the $G$-crossed graphical calculus that differ from the usual one for MTCs.

Diagrams are read from bottom to top and interpreted as quantum mechanical processes between anyons and symmetry defects. For simplicity, in what follows we assume our theory is multiplicity-free and suppress the vertex labels. Recall the notation $^ha_g$ is shorthand for $\rho_h(a_g)$, and $\bar{g}=g^{-1}$.

Trivalent vertices encode admissible splitting/fusion processes between anyons and defects, implicitly encoding the $N$-matrix entries $N_{c_{gh}}^{a_gb_g}$. 
$$
\pspicture[shift=-0.8](-0.1,-0.5)(1.5,1.4)
  \small
  \psset{linewidth=0.5pt,linecolor=black,arrowscale=1.5,arrowinset=0.15}
  \psline{->}(0.7,0)(0.7,0.45)
  \psline(0.7,0)(0.7,0.55)
  \psline(0.7,0.55) (0.25,1)
  \psline{->}(0.7,0.55)(0.3,0.95)
  \psline(0.7,0.55) (1.15,1)	
  \psline{->}(0.7,0.55)(1.1,0.95)
  \rput[bl]{0}(0.5,-0.3){$c_{ gh}$}
  \rput[br]{0}(1.4,1.05){$b_{h}$}
  \rput[bl]{0}(-0.1,1.05){$a_{g}$}
 \scriptsize
  \endpspicture
$$

The $F$-symbols have the same diagrammatic definition as for anyons, except now the topological charges carry group element labels. 

  $$\pspicture[shift=-1.0](0,-0.45)(1.8,1.8)
  \small
  \psset{linewidth=0.5pt,linecolor=black,arrowscale=1.5,arrowinset=0.15}
  \psline(0.2,1.5)(1,0.5)
  \psline(1,0.5)(1,0)
  \psline(1.8,1.5) (1,0.5)
  \psline(0.6,1) (1,1.5)
   \psline{->}(0.6,1)(0.3,1.375)
   \psline{->}(0.6,1)(0.9,1.375)
   \psline{->}(1,0.5)(1.7,1.375)
   \psline{->}(1,0.5)(0.7,0.875)
   \psline{->}(1,0)(1,0.375)
   \rput[bl]{0}(0.05,1.6){$a_g$}
   \rput[bl]{0}(0.95,1.6){$b_h$}
   \rput[bl]{0}(1.75,1.6){${c_k}$}
   \rput[bl]{0}(0.3,0.5){$e_{gh}$}
   \rput[bl]{0}(0.9,-0.5){$d_{ghk}$}
 \scriptsize
   
  \endpspicture
= \sum_{f_{hk}} \left[F_{d_{ghk}}^{a_gb_hc_k}\right]_{f_{hk},e_{gh}}
 \pspicture[shift=-1.0](0,-0.45)(1.8,1.8)
  \small
  \psset{linewidth=0.5pt,linecolor=black,arrowscale=1.5,arrowinset=0.15}
  \psline(0.2,1.5)(1,0.5)
  \psline(1,0.5)(1,0)
  \psline(1.8,1.5) (1,0.5)
  \psline(1.4,1) (1,1.5)
   \psline{->}(0.6,1)(0.3,1.375)
   \psline{->}(1.4,1)(1.1,1.375)
   \psline{->}(1,0.5)(1.7,1.375)
   \psline{->}(1,0.5)(1.3,0.875)
   \psline{->}(1,0)(1,0.375)
   \rput[bl]{0}(0.05,1.6){$a_g$}
   \rput[bl]{0}(0.95,1.6){$b_h$}
   \rput[bl]{0}(1.75,1.6){${c_k}$}
   \rput[bl]{0}(1.25,0.45){$f_{hk}$}
   \rput[bl]{0}(0.9,-0.5){$d_{ghk}$}
 \scriptsize
  \endpspicture $$

As we saw above, the algebraic data of a $G$xBFC differs in three main ways: the $R$-symbols are extended to become the matrix entries of a $G$-braiding, and new symbols $U$ and $\eta$ arise. 

The $G$-crossed $R$-symbols are written in the same way as in a UMTC, although they are no longer the matrix entries of a braiding isomorphism $c_{X,Y}:X \otimes Y \to Y \otimes X$ but a $G$-crossed braiding, $c_{X,Y}: X \otimes Y \to \,^gY \otimes X$. 

$$\pspicture[shift=-0.75](-0.1,-0.4)(1.3,1.4)
\small
  \psset{linewidth=0.5pt,linecolor=black,arrowscale=1.5,arrowinset=0.15}
  \psline(0.96,0.05)(0.2,1)
  \psline{->}(0.96,0.05)(0.28,0.9)
  \psline(0.24,0.05)(1,1)
  \psline[border=2pt]{->}(0.24,0.05)(0.92,0.9)
  \rput[bl]{0}(-0.1,1.1){$a_{g}$}
  \rput[br]{0}(1.2,1.1){$b_{h}$}
  \rput[bl]{0}(-0.15,-0.05){$b_{ h}$}
  \rput[br]{0}(1.55,-0.1){$^{ \bar h} a_{g}$}
  \endpspicture$$

Crossings of defect lines are resolved in the same manner as for anyons, except now the $R$-symbols are $G$-crossed.

$$
\pspicture[shift=-0.8](-0.1,-0.5)(1.5,1.4)
  \small
  \psset{linewidth=0.5pt,linecolor=black,arrowscale=1.5,arrowinset=0.15}
  \psline{->}(0.7,0)(0.7,0.43)
  \psline(0.7,0)(0.7,0.5)
 \psarc(0.8,0.6732051){0.2}{120}{240}
 \psarc(0.6,0.6732051){0.2}{-60}{35}
  \psline (0.6134,0.896410)(0.267,1.09641)
  \psline{->}(0.6134,0.896410)(0.35359,1.04641)
  \psline(0.7,0.846410) (1.1330,1.096410)	
  \psline{->}(0.7,0.846410)(1.04641,1.04641)
  \rput[bl]{0}(0.5,-0.3){$c_{gh}$}
  \rput[br]{0}(1.4,1.15){$b_{ h}$}
  \rput[bl]{0}(-0.1,1.15){$a_{ g}$}
 \scriptsize
  \endpspicture
=
R_{c_{gh}}^{a_{g} b_{ h}}
\pspicture[shift=-0.8](-0.1,-0.5)(1.5,1.4)
  \small
  \psset{linewidth=0.5pt,linecolor=black,arrowscale=1.5,arrowinset=0.15}
  \psline{->}(0.7,0)(0.7,0.45)
  \psline(0.7,0)(0.7,0.55)
  \psline(0.7,0.55) (0.25,1)
  \psline{->}(0.7,0.55)(0.3,0.95)
  \psline(0.7,0.55) (1.15,1)	
  \psline{->}(0.7,0.55)(1.1,0.95)
  \rput[bl]{0}(0.5,-0.3){$c_{ gh}$}
  \rput[br]{0}(1.4,1.05){$b_{ h}$}
  \rput[bl]{0}(-0.1,1.05){$a_{ g}$}
 \scriptsize
  \endpspicture
$$
The picture is similar for a left-handed crossing, but with a factor of $\left(R_{c_{gh}}^{a_gb_h}\right)^{-1}$. Now instead of satisfying a hexagon, the $G$-crossed $R$-symbols satisfy $G$-crossed heptagon equations \cite{BBCW}.
$$\includegraphics[scale=0.55]{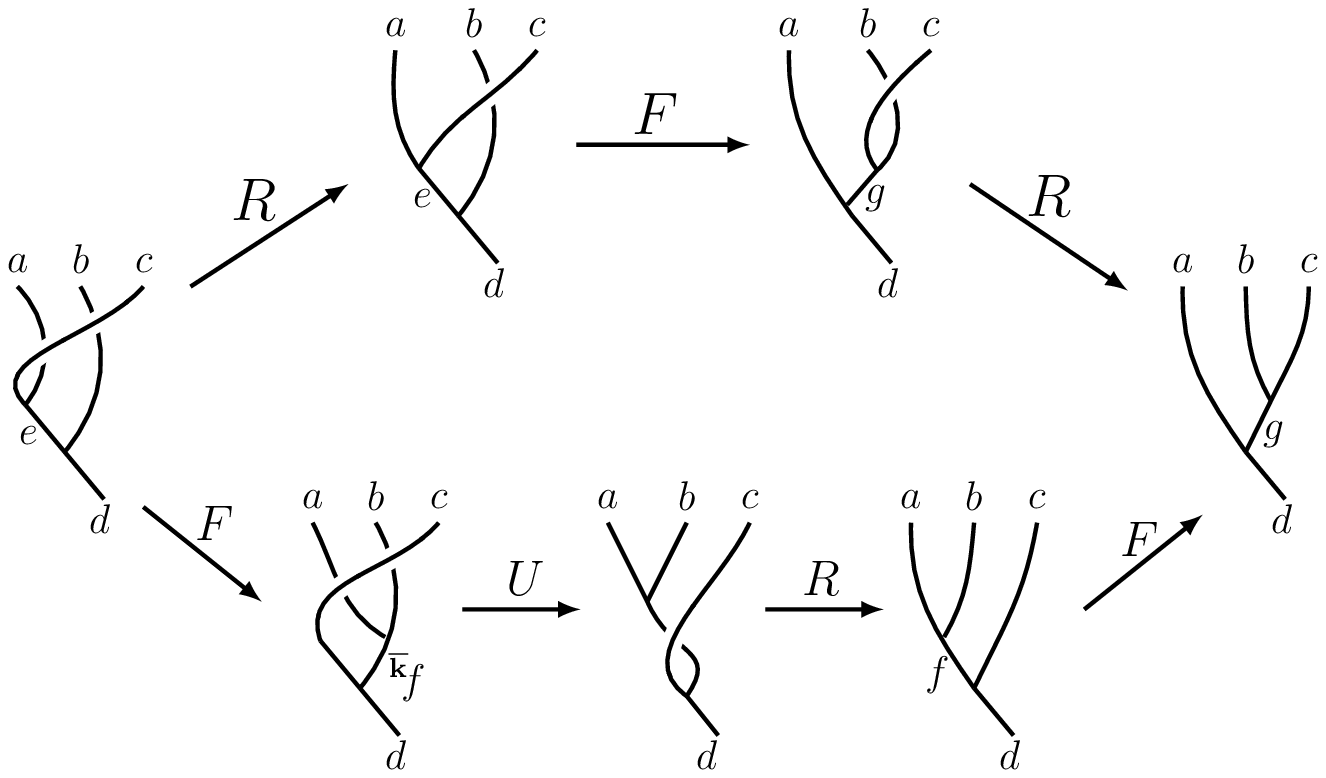}$$
A similar picture encodes compatibility of the associators with $R^{-1}$.

While for anyons strands can be freely moved over and under vertices, $U$- and $\eta$- symbols are picked up when a strand labeled by a symmetry defect passes over or under a trivalent vertex, respectively. Sliding a defect over a splitting vertex results in a $U$-symbol as follows.

$$
\label{eq:GcrossedU}
\psscalebox{1}{
\pspicture[shift=-1.7](-0.8,-0.8)(1.8,2.4)
  \small
  \psset{linewidth=0.5pt,linecolor=black,arrowscale=1.5,arrowinset=0.15}
  \psline{->}(0.7,0)(0.7,0.45)
  \psline(0.7,0)(0.7,0.55)
  \psline(0.7,0.55)(0.25,1)
  \psline(0.7,0.55)(1.15,1)	
  \psline(0.25,1)(0.25,2)
  \psline{->}(0.25,1)(0.25,1.9)
  \psline(1.15,1)(1.15,2)
  \psline{->}(1.15,1)(1.15,1.9)
  \psline[border=2pt](-0.65,0)(2.05,2)
  \psline{->}(-0.65,0)(0.025,0.5)
  \rput[bl]{0}(-0.4,0.6){$x_{ k}$}
  \rput[br]{0}(1.5,0.65){$^{\bar k}b$}
  \rput[bl]{0}(0.4,-0.4){$^{\bar k}c_{ gh}$}
  \rput[br]{0}(1.3,2.1){$b_{ h}$}
  \rput[bl]{0}(0.0,2.05){$a_{g}$}
 \scriptsize
  \endpspicture
}
=
U_{ k}\left(a, b; c\right)
\psscalebox{1}{
\pspicture[shift=-1.7](-0.8,-0.8)(1.8,2.4)
  \small
  \psset{linewidth=0.5pt,linecolor=black,arrowscale=1.5,arrowinset=0.15}
  \psline{->}(0.7,0)(0.7,0.45)
  \psline(0.7,0)(0.7,1.55)
  \psline(0.7,1.55)(0.25,2)
  \psline{->}(0.7,1.55)(0.3,1.95)
  \psline(0.7,1.55) (1.15,2)	
  \psline{->}(0.7,1.55)(1.1,1.95)
  \psline[border=2pt](-0.65,0)(2.05,2)
  \psline{->}(-0.65,0)(0.025,0.5)
  \rput[bl]{0}(-0.4,0.6){$x_{ k}$}
  \rput[bl]{0}(0.4,-0.4){$^{ \bar k}c_{gh}$}
  \rput[bl]{0}(0.15,1.2){$c_{ gh}$}
  \rput[br]{0}(1.3,2.1){$b_{h}$}
  \rput[bl]{0}(0.0,2.05){$a_{ g}$}
 \scriptsize
  \endpspicture
  }$$

$$
\psscalebox{1}{
\pspicture[shift=-1.7](-0.8,-0.8)(1.8,2.4)
  \small
  \psset{linewidth=0.5pt,linecolor=black,arrowscale=1.5,arrowinset=0.15}
  \psline{->}(0.7,0)(0.7,0.45)
  \psline(0.7,0)(0.7,1.55)
  \psline(0.7,1.55)(0.25,2)
  \psline{->}(0.7,1.55)(0.3,1.95)
  \psline(0.7,1.55) (1.15,2)	
  \psline{->}(0.7,1.55)(1.1,1.95)
  \psline[border=2pt](-0.65,2)(2.05,0)
  \psline{->}(2.05,0)(1.375,0.5)
  \rput[bl]{0}(1.4,0.6){$x_{ k}$}
  \rput[bl]{0}(0.4,-0.4){$^{ k}c_{ gh}$}
  \rput[bl]{0}(0.85,1.2){$c_{gh}$}
  \rput[br]{0}(1.3,2.1){$b_{h}$}
  \rput[bl]{0}(0.0,2.05){$a_{ g}$}
 \scriptsize
  \endpspicture
}
=
U_{ k}\left(\,^{ k}a,\,^{ k}b;\,^{k}c\right)
\psscalebox{1}{
\pspicture[shift=-1.7](-0.8,-0.8)(1.8,2.4)
  \small
  \psset{linewidth=0.5pt,linecolor=black,arrowscale=1.5,arrowinset=0.15}
  \psline{->}(0.7,0)(0.7,0.45)
  \psline(0.7,0)(0.7,0.55)
  \psline(0.7,0.55)(0.25,1)
  \psline(0.7,0.55)(1.15,1)	
  \psline(0.25,1)(0.25,2)
  \psline{->}(0.25,1)(0.25,1.9)
  \psline(1.15,1)(1.15,2)
  \psline{->}(1.15,1)(1.15,1.9)
  \psline[border=2pt](-0.65,2)(2.05,0)
  \psline{->}(2.05,0)(1.375,0.5)
  \rput[bl]{0}(1.4,0.6){$x_{ k}$}
  \rput[bl]{0}(0.4,-0.4){$^{k}c_{ gh}$}
  \rput[br]{0}(1.3,2.1){$b_{ h}$}
  \rput[bl]{0}(0.0,2.05){$a_{g}$}
  \rput[bl]{0}(-0.2,0.8){$^{ k}a$}
  \rput[bl]{0}(0.8,0.3){$^{ k}b$}
 \scriptsize
  \endpspicture
}
\qquad
$$
The picture is similar for fusion vertices, see \cite{BBCW}. 
Defect lines can be slid under splitting vertices at the cost of an $\eta$-symbol as follows.
$$\psscalebox{1}{
\pspicture[shift=-1.7](-0.8,-0.8)(1.8,2.4)
  \small
  \psset{linewidth=0.5pt,linecolor=black,arrowscale=1.5,arrowinset=0.15}
  \psline(-0.65,0)(2.05,2)
  \psline[border=2pt](0.7,0.55)(0.25,1)
  \psline[border=2pt](1.15,1)(1.15,2)
  \psline(0.7,0.55)(1.15,1)	
  \psline{->}(0.7,0)(0.7,0.45)
  \psline(0.7,0)(0.7,0.55)
  \psline(0.25,1)(0.25,2)
  \psline{->}(0.25,1)(0.25,1.9)
  \psline{->}(1.15,1)(1.15,1.9)
  \psline{->}(-0.65,0)(0.025,0.5)
  \rput[bl]{0}(-0.5,0.6){$x_{ k}$}
  \rput[bl]{0}(0.5,1.3){$^{\bar g}x$}
  \rput[bl]{0}(1.6,2.1){$^{ \overline{gh}}x_{ k}$}
  \rput[bl]{0}(0.4,-0.35){$c_{ gh}$}
  \rput[br]{0}(1.3,2.1){$b_{ h}$}
  \rput[bl]{0}(0.0,2.05){$a_{ g}$}
 \scriptsize
  \endpspicture
}
=
\eta_{x}\left({ g},{ h}\right)
\psscalebox{1}{
\pspicture[shift=-1.7](-0.8,-0.8)(1.8,2.4)
  \small
  \psset{linewidth=0.5pt,linecolor=black,arrowscale=1.5,arrowinset=0.15}
  \psline(-0.65,0)(2.05,2)
  \psline[border=2pt](0.7,0)(0.7,1.55)
  \psline{->}(0.7,0)(0.7,0.45)
  \psline(0.7,1.55)(0.25,2)
  \psline{->}(0.7,1.55)(0.3,1.95)
  \psline(0.7,1.55) (1.15,2)	
  \psline{->}(0.7,1.55)(1.1,1.95)
  \psline{->}(-0.65,0)(0.025,0.5)
  \rput[bl]{0}(-0.5,0.6){$x_{ k}$}
  \rput[bl]{0}(1.6,2.1){$^{\overline{gh}}x_{ k}$}
  \rput[bl]{0}(0.4,-0.35){$c_{ gh}$}
  \rput[br]{0}(1.3,2.1){$b_{ h}$}
  \rput[bl]{0}(0.0,2.05){$a_{ g}$}
 \scriptsize
  \endpspicture
}
$$
$$
\psscalebox{1}{
\pspicture[shift=-1.7](-0.8,-0.8)(1.8,2.4)
  \small
  \psset{linewidth=0.5pt,linecolor=black,arrowscale=1.5,arrowinset=0.15}
  \psline(2.05,0)(-0.65,2)
  \psline[border=2pt](0.7,0)(0.7,1.55)
  \psline{->}(0.7,0)(0.7,0.45)
  \psline(0.7,1.55)(0.25,2)
  \psline{->}(0.7,1.55)(0.3,1.95)
  \psline(0.7,1.55) (1.15,2)	
  \psline{->}(0.7,1.55)(1.1,1.95)
  \psline{->}(2.05,0)(1.375,0.5)
  \rput[bl]{0}(-0.8,2.1){$x_{ k}$}
  \rput[bl]{0}(1.4,0.6){$^{\overline{gh}} x_{ k}$}
  \rput[bl]{0}(0.4,-0.35){$c_{ gh}$}
  \rput[br]{0}(1.3,2.1){$b_{ h}$}
  \rput[bl]{0}(0.0,2.05){$a_{ g}$}
 \scriptsize
  \endpspicture
}
=
\eta_{x}\left({ g},{ h}\right)
\psscalebox{1}{
\pspicture[shift=-1.7](-0.8,-0.8)(1.8,2.4)
  \small
  \psset{linewidth=0.5pt,linecolor=black,arrowscale=1.5,arrowinset=0.15}
  \psline(2.05,0)(-0.65,2)
  \psline[border=2pt](0.7,0.55)(1.15,1)
  \psline[border=2pt](0.25,1)(0.25,2)
  \psline[border=2pt](1.15,1)(1.15,2)
  \psline(0.7,0.55)(0.25,1)
  \psline{->}(0.7,0)(0.7,0.45)
  \psline(0.7,0)(0.7,0.55)
  \psline(0.25,1)(0.25,2)
  \psline{->}(0.25,1)(0.25,1.9)
  \psline{->}(1.15,1)(1.15,1.9)
  \psline{->}(2.05,0)(1.375,0.5)
  \rput[bl]{0}(-0.8,2.1){$x_{ k}$}
  \rput[bl]{0}(1.4,0.6){$^{\overline{gh}} x_{ k}$}
  \rput[bl]{0}(0.5,1.3){$^{\bar g}x$}
  \rput[bl]{0}(0.4,-0.35){$c_{ gh}$}
  \rput[br]{0}(1.3,2.1){$b_{ h}$}
  \rput[bl]{0}(0.0,2.05){$a_{g}$}
 \scriptsize
  \endpspicture
}
$$
As with the $U$-symbols there is a similar picture for sliding defect lines under fusion vertices. Notice that the orientation of the $g$-defect with respect to the the trivalent vertex results in different arguments to the $U$- and $\eta$-symbols. 

In practice, understanding these moves of diagrams is sufficient to work in the $G$-extended graphical calculus. The $G$-crossed consistency equations enforcing compatibility of all of these structures are quite complicated, and the full set of such equations is listed in \cite{BBCW}, although there is some redundancy. In particular, many of the diagrams involving compatibility of the $U$- and $\eta$- symbols can be derived from the $G$-crossed pentagons and heptagons, and hence the consistency equations do not encode additional information \cite{BT}.

\section{Topological quantum computing with symmetry defects}

The mathematical formalism of topological quantum computation using UGxBFCs proceeds along the same lines as for UMTCs. Qudits are encoded in the same manner, namely as subspaces of $\Hom$ spaces of the category, and the associated braid group representations have essentially the same definition, only now they are projective. In this section we describe how to encode quantum information in the state spaces of collections of symmetry defects and how to produce quantum gates from braiding symmetry defects. To illustrate the theory we analyze some quantum gate sets arising from the examples of Section 2.

\subsection{Qudit encodings}
Let $\CC^{\times}_G=\bigoplus_{g \in G} \CC_g$ denote a UGxBFC. Let $x_{g_i}$ denote $n+1$ simple objects in $\CC_{g_i}$, with $0 \le i \le n$. Then $V_{x_{g_0}, \ldots, x_{g_n}}=\Hom\left(x_{g_0}, \bigotimes_{i=1}^n x_{g_i}  \right)$ is a finite-dimensional vector space that represents the state space of the objects (anyons and defects) $\{x_{g_i}\}_{i \ge 1}$ with total charge $x_{g_0}$. Typically for applications to quantum computing one takes all topological charges to be identical, and considers vector spaces of the form $V_{n,x_g}^{y_h}:=\Hom(y_h, x_g^{\otimes n})$.

Then $V_{n,x_g}^{y_h}$ is the Hilbert space of states associated to a collection of $n$ symmetry defects $x_g \in \CC_g$ with total charge $y_h \in \CC_h$ on a surface with trivial topology. We identify a $d$-dimensional qudit $\mathds{C}^{d}$ with $V_{n,x_g}^{y_h}$. 
\begin{figure}[H] 
 \begin{center}
$V_{n,x_g}^{y_h}= \hspace{2pt} $\mymk{$x_g \hspace{5pt} x_g \hspace {10pt} \cdots \hspace{10pt} x_g$ }\vspace{-5pt}\hspace{5pt}$y_h$
\end{center}
\end{figure}
As we will see in the next subsection these vector spaces give rise to a projective representation of the $n$-strand braid group $\BB_n$. 

$$\BB_n = \Big \langle \sigma_1, \sigma_2, \ldots, \sigma_{n-1} \Bigm\lvert \begin{array}{l} \sigma_i\sigma_j = \sigma_j \sigma_i \text{ for } |i -j| \ge 2 \\ \sigma_i\sigma_{i+1}\sigma_i = \sigma_{i+1}\sigma_i\sigma_{i+1}, i=1, 2, \ldots, n-1 \end{array}\Big \rangle.$$

In general $V_{x_{g_0}, \ldots, x_{g_n}}$ will not be a representation of the braid group, or even the pure braid group. However, it is a representation of the subgroup of $\BB_n$ that preserves the data of the topological charges $x_{g_i}$.

As with anyons, for physical reasons one is typically interested in state spaces with trivial total charge. When all $n$ simple objects $x_{g_i}$ live in $\CC_0$, i.e.\ are anyons, then the formalism reduces to TQC with anyons, see \cite{CBMS}.

\subsection{Projective braid group representations}
Next we define the representations associated to $n$ symmetry defects $x_g$ with total charge $y_h$. When the $n$ topological charges are identical, there is an action of the $n$-strand braid group on $\Hom(y_h, x_g^{\otimes n})$.
An orthonormal basis of qudit states is given by fusion trees whose internal edges are admissibly labeled by simple objects in $\CC^{\times}_G$, normalized with respect to the Markov trace inner product.  

$$\begin{tikzpicture}[scale=.25,baseline=-10, line width=.5]

\draw (1,-3/2)--(8,-5);
\draw (3,-3/2)--(3,-5/2);
\draw (4.5,-2) node {$\cdots$};
\draw (6,-3/2)--(6,-8/2);
\draw (8,-3/2)--(8, -6);
\draw (1,-3/2) node[above] {$x_g$};
\draw (3,-3/2) node[above] {$x_g$};
\draw (6,-3/2) node[above] {$x_g$};
\draw (8,-3/2) node[above] {$x_g$};

\draw (8,-6.5) node[below] {$y_h$};
\end{tikzpicture}$$
In the remainder of the paper we suppress the arrows indicating the orientations when drawing fusion trees, and assume all diagrams are oriented from the bottom up.  

The action of the braid group is defined exactly the same as for anyons. Abstractly, each element of the braid group is interpreted as a morphism $b: x_g^{\otimes n}\to x_g^{\otimes n}$. Then elements of the braid group act on a fusion tree basis element $f \in V_{n,x_g}^{y_h}$ by postcomposition:

$$b\cdot f = b\circ f.$$

On the level of diagrams, elements of the braid group $b \in \BB_n$ act on an admissibly labeled fusion trees by stacking the diagram $\begin{tikzpicture}[scale=.35,baseline=0] \draw(0,0) rectangle (1,1); \draw (1/2,1/2) node {$b$}; \end{tikzpicture}$ of $b$ on the fusion tree. 

$$\begin{tikzpicture}[scale=.25,baseline=-10, line width=.5]
\draw (4.5,2) node {$b$};
\draw (1,-3/2) rectangle (8,11/2);
\draw (1,-3/2)--(8,-5);
\draw (3,-3/2)--(3,-5/2);
\draw (4.5,-2) node {$\cdots$};
\draw (6,-3/2)--(6,-8/2);
\draw (8,-3/2)--(8, -6);
\draw (1,-1) node {};
\draw (2,-1) node  {};
\draw (3,-1) node  {};
\draw (4,-1) node  {};
\draw (5,-1) node  {};
\draw (6,-1) node  {};
\draw (7,-1) node  {};
\draw (8,-1) node  {};
\draw (2.25, -2.75) node  {};
\draw (3.5, -3.25) node {};
\draw(4.5,-3.70) node  {};
\draw(5.5,-4.25) node  {};
\draw(6.5,-4.75) node {};
\draw(7.5,-5.25) node {};
\draw (8,-6.5) node {};
\end{tikzpicture}$$
By extending linearly to all of $V_{n,x_g}^{y_h}$ this gives a projective representation $\rho_{n, x_g, y_h}$, which we refer to as the representation corresponding to $n$ symmetry defects of type $x_g$ with total charge $y_h$. 

The matrix representation in the fusion tree basis is then given by resolving $b\cdot f$ as a linear combination of admissibly labeled fusion trees using the graphical calculus. (Or in practice, by taking inner products as in Section 4.) However, this defines a braid group representation which is not linear in general, only projective. This is due to symmetry defects being confined, see \cite{BBCW}. However, since quantum mechanics is insensitive to a global phase already, even when analyzing representations coming from UMTCs one is only interested in the projective image. That means we can meaningfully compare quantum computers designed using either anyons or symmetry defects.

\subsubsection{Universality}
For applications to quantum computing one desires to be able to approximate any unitary matrix in $U(d)$ to arbitrary accuracy with only polynomially many matrices in some small gate set $S \subset U(d)$. 

For gate sets coming from braiding $n$ symmetry defects $x_g$ with total charge $y_h$, this translates into the question of whether the projective image of the braid group representation $V_{n,x_g}^{y_h}$ is dense in the unitary group. That is, one wants $\overline{{\rho}_{n, x_g, y_h}(\BB_n)} \supset SU(d)$ projectively. 

Thus the problem of analyzing the computational power of a symmetry defect model boils down to studying the images of its associated projective braid group representations in the unitary group. As a result all of the techniques and results used to study quantum gates in general quantum information apply. In particular, standard results about optimality of encodings and leakage out of computational subspaces carry over without modification \cite{AS}.

For application to quantum information one also needs to know how to represent non-unitary physical processes, for example projective measurement of topological charge.

\subsubsection{Projective measurement}
Not only is projective measurement the way to perform read-out of a topological quantum computer, but it is also used in the bulk of computation as a way to force measurement outcomes of subroutines. Projective measurement can also be leveraged to simulate physical braiding along with the help of ancillary qudits \cite{BFN}. This may be especially relevant for symmetry defects, which behave differently energetically than anyons and whose physical realizations may not be as conducive to physical braiding \cite{BBCW,BJQ}. As it stands many state-of-the-art proposals for anyonic hardware rely on the measurement-only formalism, see for example \cite{KK}.

On the level of diagrams, a projective measurement of symmetry defects $a_g$ and $b_h$ onto $c_{gh}$ is given by 

$$\pspicture[shift=-0.6](-0.1,-0.45)(1.4,1)
  \small
  \psset{linewidth=0.5pt,linecolor=black,arrowscale=1.5,arrowinset=0.15}
  \psline{->}(0.7,0)(0.7,0.45)
  \psline(0.7,0)(0.7,0.55)
  \psline(0.7,0.55) (0.25,1)
  \psline{->}(0.7,0.55)(0.3,0.95)
  \psline(0.7,0.55) (1.15,1)
  \psline{->}(0.7,0.55)(1.1,0.95)
  \rput[bl]{0}(0.2,0.2){$c_{gh}$}
  \rput[br]{0}(1.5,0.8){$b_h$}
  \rput[bl]{0}(-.1,0.8){$a_g$}
  \psline(0.7,0) (0.25,-0.45)
  \psline{-<}(0.7,0)(0.35,-0.35)
  \psline(0.7,0) (1.15,-0.45)
  \psline{-<}(0.7,0)(1.05,-0.35)
  \rput[br]{0}(1.5,-0.5){$b_h$}
  \rput[bl]{0}(-.1,-0.5){$a_g$}
\scriptsize
  \endpspicture.$$
  
 The matrix representation of this morphism can be calculated by specifying a basis of fusion trees, then stacking the diagram and resolving it as a linear combination of basis elements via the rules of the graphical calculus, exactly as one would for a unitary process. The example in Section 4 incorporates this technique.

Having covered the construction of quantum gates from symmetry defects we return to the small examples introduced in the previous section.

\subsection{Toric code with electromagnetic duality}

The $F$-, $G$-crossed $R$-, $U$-, and $\eta$-symbols were found in \cite{BBCW} by making a gauge choice and solving the $G$-crossed consistency equations. It was also observed that the fusion subcategories $\{1, \psi, \sigma_+\}$ and $\{1, \psi, \sigma_-\}$ have Ising fusion rules. We show that this relationship extends even further: a toric code defect qubit is identical to an Ising qubit.  

We consider the projective representation of the four-strand braid group $\BB_4$ afforded by braiding four symmetry defects $\sigma_+$ with trivial total charge. The fusion space $V_{4,\sigma_+}^1=\Hom(1, \sigma_+^{\otimes 4})$ is two-dimensional, corresponding to the admissible labelings of the following fusion tree.

$$|a \rangle = \begin{tikzpicture}[yshift=5mm, baseline=-45,line width=.5]
\draw (1,-1.5) node[above] {$\sigma_+$}--(2.5,-3)--(4,-1.5) node[above] {$\sigma_+$};
\draw (2,-1.5) node[above] {$\sigma_+$}--(1.5,-2);
\draw (3,-1.5) node[above] {$\sigma_+$}--(3.5,-2);
\draw (2.5,-3)--(2.5,-3.5) node[below] {$1$};
\draw (2,-2.75) node[left] {$a$};
\draw (3,-2.75) node[right] {$a$};
\end{tikzpicture}$$
Thus the symmetry defects encode a qubit with basis states $|1\rangle$ and $|\psi\rangle$. 

We find
$$\rho(\sigma_1)=\rho(\sigma_3)=e^{-\pi i/8} \begin{pmatrix} 1 & 0 \\ 0 & i \\ \end{pmatrix}$$
and
$$\rho(\sigma_2)=\frac{e^{-\pi i/8}}{2} \begin{pmatrix}1+i & i-1\\ 1-i & 1+i \\ \end{pmatrix}.$$

We remark that this generates the same image as the Ising qubit \cite{CBMS}.
With the gauge choice made here the generating gates are actually equal on the nose to the Ising qubit gates. The Ising qubit is a Jones representation of $\BB_4$, whose projective image is finite and given by $\mathds{Z}_2^2 \rtimes S_3$ \cite{J83}.

\subsection{$\left(\mathds{Z}_3\right)^X_{\mathds{Z}_2}$ symmetry defect qutrit}

We consider the projective representation of the four-strand braid group $\BB_4$ afforded by braiding four symmetry defects $X_{\omega}$ with trivial total charge. The fusion space $V_{4,X_{\omega}}^{1}=\Hom(1, X_{\omega}^4)$ is three-dimensional, corresponding to the admissible labelings of the following fusion tree.

$$|ab\rangle=\begin{tikzpicture}[yshift=5mm, baseline=-45, line width=.5]
\draw (1,-1.5) node[above] {$X_{\omega}$}--(2.5,-3)--(4,-1.5) node[above] {$X_{\omega}$};
\draw (2,-1.5) node[above] {$X_{\omega}$}--(1.5,-2);
\draw (3,-1.5) node[above] {$X_{\omega}$}--(3.5,-2);
\draw (2.5,-3)--(2.5,-3.5) node[below] {$1$};
\draw (2,-2.75) node[left] {$a$};
\draw (3,-2.75) node[right] {$b$};
\end{tikzpicture}$$

We work in the ordered basis given by $\{|11\rangle, |\omega \omega^*\rangle, |\omega^*\omega \rangle\}$. We find
$$\rho(\sigma_1) =\rho(\sigma_3)\sim \begin{pmatrix} R^{X_{\omega}X_{\omega}}_1 & 0 & 0 \\ 0 & R^{X_{\omega}X_{\omega}}_{\omega} & 0 \\ 0 & 0 & R^{X_{\omega}X_{\omega}}_{\omega^*} \end{pmatrix}$$
and
$$\rho(\sigma_2) \sim \frac{1}{3\sqrt{3}}\begin{pmatrix} 1 & \frac{1}{2}(-1+i\sqrt{3}) & \frac{1}{2}(-1+i\sqrt{3}) \\ \frac{1}{2}(-1+i\sqrt{3}) & 1 & \frac{1}{2}(-1+i\sqrt{3}) \\\frac{1}{2}(-1+i\sqrt{3}) & \frac{1}{2}(-1+i\sqrt{3}) & 1 \end{pmatrix}.$$

We remark that these matrices form a subset of the matrices generating the Jones representation of $\BB_4$ coming from anyon 1 and total charge 2 in the $\ZZ_2$-gauged theory $\ZZ_3^{\ZZ_2}=TL_4$, the Temperley-Lieb-Jones algebroid at $A=i e^{-\pi i /12}$. In particular the image of the projective representation is finite \cite{J83}.

So far all of our examples have been concerned with what happens when one is interested in defect-defect braiding alone. Next we go more in depth with an example that shows how to incorporate anyon-defect braiding and measurement into the mathematical formalism.

\section{Bilayer Ising with $\mathds{Z}_2$ layer-exchange symmetry}

In this section we outline how to get a $T$-gate from a $\mathds{Z}_2$-extension of the category $\textbf{Ising}^{(1)} \boxtimes \textbf{Ising}^{(1)}$. This follows the work of physicists who argued that certain types of topological defects called ``genons'' that arise in bilayer Ising anyon models can have non-abelian statistics \cite{BJQ}. It is important to note that there are eight inequivalent MTCs of rank 3 with nonintegral quantum dimensions which all go by the name of Ising. These are the categories $\textbf{Ising}^{(\nu)}$ with isomorphism classes of simple objects $\{1,\sigma,\psi\}$, where $\nu$ parametrizes the quantum twist of the anyon $\sigma$. These lead to 20 inequivalent MTCs of the form $\textbf{Ising}^{(\nu_1)} \boxtimes \textbf{Ising}^{(\nu_2)}$ \cite{BGPR}. Here we take $\nu_1=\nu_2=1$, and by \textbf{Ising} we mean $\textbf{Ising}^{(1)}$, the category with modular data given in \cite{CBMS}.
 
 \subsection{Towards a universal quantum computer based on Ising anyons}
Prior to the description of symmetry defects using UGxBFCs \cite{BBCW}, physicists argued that certain types of symmetry defects called \emph{genons} could give rise to non-abelian braiding statistics and thus could be used together with anyons for quantum information processing \cite{BJQ}. 
 
Braiding Ising anyons alone is not universal for quantum computation: in particular, the $T$-gate is missing for universal single-qubit operations. The authors of \cite{BJQ} designed a protocol that addressed this shortcoming in another way, by constructing a $T$-gate from braiding and measurement of $\mathds{Z}_2$ symmetry defects associated to two layers of the Ising theory with layer-exchange symmetry. More precisely, they furnished a universal gate set $\{T=G_1, G_2, G_3\}$ which can be realized by braiding anyons and defects. Mathematically, this is described by a $\mathds{Z}_2$-extension of the Deligne product of the category $\textbf{Ising}$ with itself, $\left(\textbf{Ising} \boxtimes \textbf{Ising}\right)^{\times}_G$. 

It is an open question to find a full solution to the $G$-crossed consistency equations for $\left(\textbf{Ising} \boxtimes \textbf{Ising}\right)^{\times}_G$. However, under certain physically motivated assumptions about the algebraic data defining this $\mathds{Z}_2$-extension, we use the mathematical formalism of fusion categories to verify that the physical protocol of \cite{BJQ} produces a $T$-gate.

A mathematical treatment of this proposal proceeds by formulating the hypothetical data as a conjecture about a symmetry defect $X_1$ in the category $\left( \textbf{Ising} \boxtimes \textbf{Ising} \right)^{\times}_{\ZZ_2}$. For the motivation behind the conjecture see \cite{BJQ}.

\begin{conj} There exists a set of solutions to the $\mathds{Z}_2$-crossed consistency equations associated to  $\textbf{Ising}\boxtimes \textbf{Ising}$ such that
$$R^{X_1X_1}_{a\boxtimes a}=\theta_{a} \hspace{5pt} \text{ for } a \in \textbf{Ising}$$ and $$[F^{X_1X_1X_1}_{X_1}]=S_{\textbf{Ising}}= \frac{1}{2}\begin{pmatrix} 1 & \sqrt{2} & 1 \\ \sqrt{2} & 0 & -\sqrt{2} \\ 1 & -\sqrt{2} & 1 \\ \end{pmatrix}.$$
\end{conj}
Our analysis etablishes the following corollary of the conjecture.
\begin{cor}[To Conjecture 4.1] Braiding and projective measurement of bilayer Ising anyons together with $\mathds{Z}_2$ symmetry defects is universal for quantum computation.
\end{cor}

Should Conjecture 4.1 hold, our analysis provides mathematical justification for including symmetry defects in the quantum information scientist's toolkit for producing universal gate sets from topological phases of matter.

In the remainder of this section we outline the proof of the corollary, which involves encoding the protocol in the morphism space of the UGxBFC given by $\left( \textbf{Ising} \boxtimes \textbf{Ising} \right)^{\times}_{\ZZ_2}$ and calculating its matrix representation. Due to the projective measurement, the protocol is not a braid, and hence the $T$-gate is not in the projective braid group image. But the algebraic data of the $G$-crossed category nevertheless affords a well-defined matrix representation corresponding to the physical process.

\subsection{The category $\left( \textbf{Ising} \boxtimes \textbf{Ising} \right)^{\times}_{\ZZ_2}$}

As a $\mathds{Z}_2$-graded category, 

$$\left(\textbf{Ising}^{\boxtimes 2}\right)^{\times}_{\mathds{Z}_2}= \bigoplus_{a,b \in \{1,\sigma,\psi \}}\{a\boxtimes b\} \bigoplus \{X_1, X_{\sigma},X_{\psi}\}.$$ 

That is, there are 12 isomorphism classes of simple objects: 9 anyon types and 3 symmetry defect types. The defects are labeled $X_a$ for the fixed points $a\boxtimes a$ of the layer-exchange symmetry.

The quantum dimensions are given by
$$ d_{11}=d_{1\psi}=d_{\psi 1}=d_{\psi\psi}=1,\; d_{1 \sigma}=d_{\sigma 1}=d_{\sigma \psi}=d_{\psi \sigma}=\sqrt{2},\; d_{\sigma\sigma}=2$$
and
$$d_{X_1}=d_{X_{\psi}}=2,\;d_{X_{\sigma}}=2\sqrt{2}.$$

Fusion in the trivial sector is given by \textbf{Ising} fusion in each factor of the Deligne product. Among the defects the fusion channels are by given by
$$\begin{array}{c||c|c|c}

\otimes & X_1 & X_{\sigma} & X_{\psi} \\
\hline
\hline
 X_1 & \text{11 $\oplus $ $\sigma \sigma $ $\oplus $ $\psi \psi $} & \text{1$\sigma $ $\oplus $ $\sigma $1 $\oplus $ $\sigma \psi $ $\oplus $ $\psi \sigma $} & \text{1$\psi $ $\oplus $ $\sigma \sigma $ $\oplus $ $\psi $1} \\
 \hline
 X_{\sigma} & \text{1$\sigma $ $\oplus $ $\sigma $1 $\oplus $ $\sigma \psi $ $\oplus $ $\psi \sigma $} & \text{11 $\oplus $ 1$\psi $ $\oplus $ 2 $\sigma \sigma $ $\oplus $ $\psi $1 $\oplus $ $\psi \psi $} & \text{1$\sigma $ $\oplus $ $\sigma $1 $\oplus $ $\sigma \psi $ $\oplus $ $\psi \sigma $} \\
 \hline
 X_{\psi} & \text{1$\psi $ $\oplus $ $\sigma \sigma $ $\oplus $ $\psi $1} & \text{1$\sigma $ $\oplus $ $\sigma $1 $\oplus $ $\sigma \psi $ $\oplus $ $\psi \sigma $} & \text{11 $\oplus $ $\sigma \sigma $ $\oplus $ $\psi \psi $} \\

\end{array}$$
while between the bilayer anyons and defects they are given by

$$\begin{array}{c||c|c|c|c|c|c|c|c|c}

\otimes & 11 & 1\sigma & 1\psi & \sigma 1 & \sigma\sigma & \sigma \psi & \psi 1 & \psi \sigma & \psi \psi \\
\hline
\hline
 X_1 & X_1 & X_{\sigma } & X_{\psi } & X_{\sigma } & X_{\psi }\oplus X_1 & X_{\sigma } & X_{\psi } & X_{\sigma } & X_1 \\
 \hline
 X_{\sigma} & X_{\sigma } & X_{\psi }\oplus X_1 & X_{\sigma } & X_{\psi }\oplus X_1 & 2 X_{\sigma } & X_{\psi }\oplus X_1 & X_{\sigma } & X_{\psi }\oplus X_1 & X_{\sigma } \\
 \hline
X_{\psi} & X_{\psi } & X_{\sigma } & X_1 & X_{\sigma } & X_{\psi }\oplus X_1 & X_{\sigma } & X_1 & X_{\sigma } & X_{\psi } \\

\end{array}.$$
In what follows we will focus on the defect $X_1$ and the fusion channel
$$X_1 \otimes X_1 = 11 \oplus \sigma\sigma \oplus \psi\psi.$$
\subsection{Qubit and $T$-gate encoding}
The single-qubit $T$-gate will be realized as the logical operation enacted on a two-dimensional subspace $\mathcal{H}_{log}$ of a physical Hilbert space $\mathcal{H}_{phys}$ associated to the fusion space of a collection of four monolayer Ising anyons $\sigma 1$ and four symmetry defects $X_1$ with trivial total charge. That is, $$\mathcal{H}_{phys} = \Hom \left ( 11, \sigma1^{\otimes 4} \otimes X_1^{\otimes 4} \right)$$ as a morphism space in the category $\left( \textbf{Ising} \boxtimes \textbf{Ising} \right)^X_{\ZZ_2}$. 

We fix a basis $|x y \rangle$ of the full fusion space to be the set of admissibly labeled fusion trees, where $x \in \{11, \psi 1 \}$ and $y \in \{11, \sigma\sigma, \psi\psi\}$. 

\begin{figure}[h!]

$|xy\rangle=$
\begin{tikzpicture}[scale=.5,baseline=-60,line width=.5]

\draw (1,-3/2)--(8,-5);
\draw (2,-3/2)--(2,-2);
\draw (3,-3/2)--(3,-5/2);
\draw (4,-3/2)--(4,-3);
\draw (5,-3/2)--(5,-7/2);
\draw (6,-3/2)--(6,-8/2);
\draw (7,-3/2)--(7,-9/2);
\draw (8,-3/2)--(8, -6);
\draw (1,-1) node {$\small \sigma 1$};
\draw (2,-1) node {$\small \sigma 1$};
\draw (3,-1) node {$\small \sigma 1$};
\draw (4,-1) node {$\small \sigma 1$};
\draw (5,-1) node {$\small X_1$};
\draw (6,-1) node {$\small X_1$};
\draw (7,-1) node {$\small X_1$};
\draw (8,-1) node {$\small X_1$};
\draw (2.25, -2.75) node {\small $x$};
\draw (3.5, -3.25) node {\small $\sigma 1$};
\draw(4.5,-3.70) node {\small $11$};
\draw(5.5,-4.25) node {\small $X_1$};
\draw(6.5,-4.75) node {\small $y$};
\draw(7.5,-5.25) node {\small $X_1$};
\draw (8,-6.5) node {\small $11$};
\end{tikzpicture}
\end{figure}

Following \cite{BJQ} the logical qubit is encoded in the $\mathds{C}$-span of $\{|11, 11\rangle, |\psi 1, 11\rangle \}$. In general there will be leakage out of the logical subspace. In what follows we simplify the notation and suppress both the index label of the bottom layer and symmetry defect fusion channel and write $|11,11\rangle=|1\rangle, |\psi 1, 11 \rangle = |\psi \rangle$.
$$\mathcal{H}_{log}=\Span_{\mathds{C}} \{ |1\rangle, |\psi \rangle \}.$$

The physical protocol to initialize into this logical subspace is as follows:
\begin{enumerate}
\item  Create two pairs of monolayer Ising anyons $\sigma 1$ from the vacuum. 
\item Create two pairs of symmetry defects $X_1$ from the vacuum, using projective measurement if necessary to fix the total charge of each pair to be trivial. \\ 

The protocol to enact the logical $T$-gate is given by the following sequence of steps.\\ 
\item Select one of the monolayer Ising anyons from a pair and braid it around the middle two genons in the counterclockwise direction.
\item Fix the pair of Ising anyons from (3) to have trivial total charge with a projective measurement.
\item Perform a full exchange of the middle two defects.
\item Braid the anyon from (3) with the middle two defects in the clockwise direction.\\

Readout of the computation then proceeds by pair annihilation and projective measurement. \\

\item Annihilate the pairs of Ising anyons and the pairs of defects. 
\end{enumerate}

Steps (3)-(6) determine a morphism in $\Hom \left ( \sigma^{\otimes 4} \otimes X_1^{\otimes 4}, \sigma^{\otimes 4} \otimes X_1^{\otimes 4} \right)$. In the graphical calculus the diagram for this morphism takes the form depicted in Figure 4.1.

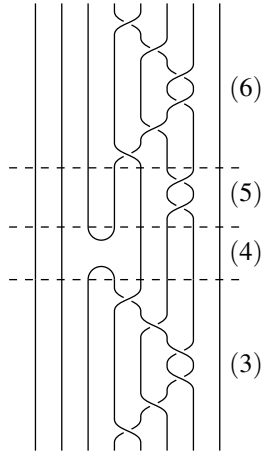
\begin{figure}[H]

$$\begin{tikzpicture}[scale=.35, line width=.5]
\braid[number of strands=8] a_4^{-1} a_5 a_6 a_6 a_5 a_4;
\draw (4,0) arc (0:180:1/2);
\draw (3,2) arc (180:360:1/2);
\draw (1,0)--(1,2);
\draw (2,0)--(2,2);
\draw (5,0)--(5,2);
\draw (6,0)--(6,2);
\draw (7,0)--(7,2);
\draw (8,0)--(8,2);

\draw[dashed] (0,0)--(9,0);
\draw[dashed] (0,2)--(9,2);
\draw[dashed] (0,4.25)--(9,4.25);

\draw (9,7.25) node {$(6)$};
\draw (9,3.25) node {$ (5)$};
\draw (9,1) node {$ (4)$};
\draw (9,-3.25) node {$(3)$};
\begin{scope}[shift={(0,10.5)}]
\braid[number of strands=8] a_4^{-1} a_5^{-1} a_6^{-1} a_6^{-1} a_5^{-1} a_4 a_6^{-1} a_6^{-1};
\end{scope}

\end{tikzpicture}$$
\caption{Diagrammatic representation of $T$-gate protocol.}
\end{figure}
\hspace{10pt}
The matrix entries $T_{ij}$ of the protocol $T$ restricted to the logical subspace can be found by calculating $\langle v_i |T | v_j \rangle$, which gives a morphism $11 \to 11$ and hence by Schur's Lemma is a scalar times the identity. 

\begin{figure}[H]
\begin{tikzpicture}[scale=.35,baseline=-10, line width=.5]

\draw (8,10) arc (180:0:1/2);
\draw (9,10)--(9,-6) arc (0:-180:1/2);
\draw (1,11/2)--(8,9);
\draw (2,11/2)--(2,6);
\draw (2.25, 6.75) node  {\tiny $1/\psi$};
\draw (3,11/2)--(3,13/2);
\draw (4,11/2)--(4,7);
\draw (5,11/2)--(5,15/2);
\draw (6,11/2)--(6,16/2);
\draw (7,11/2)--(7,17/2);
\draw (8, 11/2)--(8, 20/2);

\draw (4.5,2) node {$T$};
\draw (1,-3/2) rectangle (8,11/2);
\draw (1,-3/2)--(8,-5);
\draw (2,-3/2)--(2,-2);

\draw (3,-3/2)--(3,-5/2);
\draw (4,-3/2)--(4,-3);
\draw (5,-3/2)--(5,-7/2);
\draw (6,-3/2)--(6,-8/2);
\draw (7,-3/2)--(7,-9/2);
\draw (8,-3/2)--(8, -6);
\draw (1,-1) node {};
\draw (2,-1) node  {};
\draw (3,-1) node  {};
\draw (4,-1) node  {};
\draw (5,-1) node  {};
\draw (6,-1) node  {};
\draw (7,-1) node  {};
\draw (8,-1) node  {};
\draw (2.25, -2.75) node  {\tiny $1/\psi$};
\draw (3.5, -3.25) node {};
\draw(4.5,-3.70) node  {};
\draw(5.5,-4.25) node  {};
\draw(6.5,-4.75) node {};
\draw(7.5,-5.25) node {};
\draw (8,-6.5) node {};
\end{tikzpicture}
\caption{Diagrammatic representation of matrix entries.}
\end{figure}
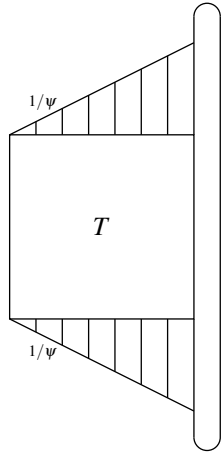
The matrix entries are determined by using the graphical calculus to find these scalars, yielding $$\rho(T) = \begin{pmatrix} \langle 1 |T| 1 \rangle & \langle \psi | T | 1 \rangle \\ \langle 1| T |\psi \rangle & \langle \psi | T| \psi \rangle \end{pmatrix}$$
up to a normalization. 

We briefly recall the tools needed for these calculations.

\subsubsection{Preliminaries for calculating matrix entries}
In addition to the $S$-matrix, $R$- and $F$- symbols for the UMTC $\textbf{Ising} \boxtimes \textbf{Ising}$, one needs (at minimum) the $R$- and $F$-symbols specified in Conjecture 4.1. Of course there are many choices of sequences of diagrammatic moves to resolve diagrams in the graphical calculus and in general one needs all of the $G$-crossed data: $\{R,F,U,\eta\}$. We note however that for this specific example we were able to calculate the matrix of $T$ in a way that is independent of $U$- and $\eta$- symbols. 

The formulas for the matrix entries involve the $S$-matrix, which is defined diagrammatically by
$$S_{ab} 
=\frac{1}{\mathcal{D}}
\pspicture[shift=-0.4](0.0,0.2)(2.6,1.3)
\small
  \psarc[linewidth=0.5pt,linecolor=black,arrows=<-,arrowscale=1.5,arrowinset=0.15] (1.6,0.7){0.5}{167}{373}
  \psarc[linewidth=0.5pt,linecolor=black,border=3pt,arrows=<-,arrowscale=1.5,arrowinset=0.15] (0.9,0.7){0.5}{167}{373}
  \psarc[linewidth=0.5pt,linecolor=black] (0.9,0.7){0.5}{0}{180}
  \psarc[linewidth=0.5pt,linecolor=black,border=3pt] (1.6,0.7){0.5}{45}{150}
  \psarc[linewidth=0.5pt,linecolor=black] (1.6,0.7){0.5}{0}{50}
  \psarc[linewidth=0.5pt,linecolor=black] (1.6,0.7){0.5}{145}{180}
  \rput[bl]{0}(0.1,0.45){$a$}
  \rput[bl]{0}(0.8,0.45){$b$}
  \endpspicture
$$
where $\mathcal{D}=\sqrt{\sum_{a \in \text{Irr}(\CC)} d_a^2}$.

\subsubsection{The diagonal entries of $\rho(T)$}

By observing that the leftmost two strands in Figure 2 are through strands, it follows from Schur's lemma that the off-diagonal entries of $\rho(T)$ vanish. That is, $$\langle \psi |T |1 \rangle = \langle 1 | T | \psi \rangle = 0. $$
As for the diagonal entries, they can be found in terms of the algebraic data of the $\mathds{Z}_2$-crossed category.
We find that up to a normalization
$$\langle 1 |T| 1 \rangle=  d_{\sigma 1} d_{X_1}^2\sum_c \left (R^{X_{1}X_{1}}_c\right )^2 \left|\left [F^{X_1X_1X_1}_{X_1}\right]_{c,11}\right|^2  \left( \frac{S_{\sigma 1,c}}{S_{11,c}} \right)^2$$

and

  $$\langle \psi|T|\psi \rangle=d_{\sigma 1}d_{X_1}^2 \sum_c  \left ( R^{X_1 X_1}_c\right )^2 \left | \left [ F^{X_1X_1X_1}_{X_1}\right ]_{c,11} \right |^2 \left ( 1 - \left ( \frac{S_{\psi 1,c}}{S_{11,c}} \right) \right),$$
  where $c \in \{11, \sigma\sigma, \psi\psi\}$. 
 The details of the diagrammatic calculations that produce these formulas are (1) somewhat subtle due to the additional rules of the $G$-crossed graphical calculus and (2) constrained by a desire to make as few assumptions as possible about the solutions to the $G$-crossed consistency equations, not yet being available. As a result the derivation of these formulas is a bit involved, so we have included it as an appendix.
 
Assuming that Conjecture 4.1 holds, one can check that

$$\langle \psi | T | \psi \rangle= e^{\pi i/4} \langle 1 |T | 1 \rangle.$$

This would confirm that the protocol implements a phase shift of $\pi/4$ between the two qubit states, effecting a $T$-gate 
$$\rho(T) \sim \begin{pmatrix} 1 & 0 \\ 0 & e^{\pi i/4} \\ \end{pmatrix}$$
on the logical qubit. This proves the corollary to the conjecture.

\section{Discussion}
When one is interested in topologically protected gate sets coming from braiding defects alone, one can analyze the image of the projective braid group representation associated to a UGxBFC as outlined in section 3. It is expected that braiding symmetry defects behaves like braiding anyons in several important ways. First, it is expected that the image of the braid group representations coming from braiding defects will be finite when the quantum dimensions are weakly integral (this is a generalization of the property-F conjecture to UGxBFCs). It is also expected that if the underlying anyons are not universal then the symmetry defects will not be universal either. 

The example of the bilayer Ising model with $\ZZ_2$ symmetry defects in Section 4 indicates how to analyze the logical operations that can be realized through braiding and measurement of anyons and symmetry defects in general. In particular, our calculation supports the previous work of physicists who argued that symmetry defects could be used along with measurement to enlarge deficient gate sets to universal ones \cite{BJQ}. The calculation provides further proof-of-concept for the inclusion of symmetry defects in the quantum information theorist's toolkit, along with projective measurement, magic state distillation, and other ways to supplement deficient anyonic gate sets. However, encoding physical protocols in the language of $G$-crossed extensions of UMTCs is evidently quite technical. 

The efforts to build a topological qubit and ultimately a scalable universal topological quantum computer are still underway. The majority of such efforts currently focus on measurement-based schemes involving Majorana zero modes (MZMs), Ising anyon-like objects that appear at the ends of nanowires placed on a superconductor in the presence of a magnetic field \cite{KK}. But it is within the realm of possibility that engineering a topological phase with symmetry defects could present a shorter path to a universal quantum computer.  

On this note, bilayer systems are interesting to study because any topological phase can be stacked into a bilayer phase that then possesses a $\ZZ_2$ bilayer-exchange symmetry. The best understood examples of topological phases are fractional quantum Hall liquids, so it is realistic to consider what experiments can be performed to probe symmetry defects in bilayer fractional quantum Hall systems.  Since fractional quantum Hall systems are electron systems, strictly speaking one needs a theory of $G$-extensions of super-modular categories in order to formulate the problem mathematically. One important exception is the $\nu=\frac{1}{3}$ fractional quantum Hall liquid, where the associated super-modular category splits into $SU(3)_1$ and the physical electron.  In this case the system can be studied via a bilayer $SU(3)_1$ topological order.  An interesting question is then how to use symmetry defects to generate a universal gate set for quantum computing from this abelian topological phase, which is available in almost every fractional quantum Hall lab. 

Apart from the physical motivations, symmetry defects in the context of TQC are interesting for mathematical reasons alone. Going forward, it is highly desirable to understand the precise mathematical relationship between the (projective) braid group representations afforded by symmetry defects and the braid group representations afforded by the underlying anyons. Relatedly, it is an open question to determine the precise relationship of the braid group representations arising from a UMTC $\CC$ and its gauging $(\CC^{\times}_G)^G$. An answer to this question could help resolve the Property-F conjecture \cite{NR}, one of most well-known open problems in modular tensor category theory. 

A systematic and detailed study of TQC with symmetry defects is to appear in \cite{DThesis}.

\section*{Acknowledgements}
We would like to thank Maissam Barkeshli, Parsa Bonderson, and Meng Cheng for contributing figures, Xingshan Cui for sharing Mathematica code, as well as Eric Samperton and the referees for helpful feedback on the draft. 

\newpage
\appendix
\appendixpage
In the appendices we provide a derivation of the matrix entries of the $T$-gate protocol from Section 4.2, with data and graphical formulas for reference.

$$\btikz[line width=.5,scale=.5]

\begin{scope}[shift={(5,6.25)}]
\braid a_1^{-1} a_1^{-1};
\end{scope}

\draw (4,0) to [out=90, in=-90] (7.5,2);
\draw (4,10) to [out=-90, in=90] (7.5,8);

\foreach \y in {6,7}
{
\draw[white,line width=7.5] (\y,0)--(\y,4);
\draw[white,line width=7.5] (\y,6.5)--(\y,10);
\draw (\y,0)--(\y,4);
\draw (\y,6)--(\y,10);
}

\draw[white, line width=7.5] (7.5,2) to [out=90, in=-90] (4,4) arc (0:180:.5)--(3,0);
\draw (7.5,2) to [out=90, in=-90] (4,4) arc (0:180:.5)--(3,0);

\draw[white, line width=7.5] (7.5,8) to [out=-90, in=90] (4,6) arc (0:-180:.5)--(3,10);
\draw (7.5,8) to [out=-90, in=90] (4,6) arc (0:-180:.5)--(3,10);


\draw[white,line width=7.5] (5,0)--(5,10);
\draw (5,0)--(5,10);
\foreach \x in {1,2,5,8}
{\draw (\x,0)--(\x,10);}

\foreach \x in {1,2,3,4}
{\draw (\x,.375) node [left] {\small $\sigma$};}
\foreach \y in {5,6,7,8}
{\draw (\y,.375) node [left] {\small $X$};}

\etikz$$

In Appendix A we list all the data needed to compute the matrix entries, including the conjectured $G$-crossed braiding and $F$-symbols for the $\mathds{Z}_2$-extension of $\bf{Ising}\boxtimes \bf{Ising}$. In Appendix B we provide a list of helpful moves in the graphical calculus for UGxBFCs (and the UMTCs they contain) for reference. In Appendix C, we demonstrate a sequence of moves in the graphical calculus that produces the formulas for the matrix entries. We find that the phase difference induced on the basis states is indeed $e^{\pi i/4}$ as anticipated by the physical argument given in \cite{BJQ}.

\section{Data}

\subsection{Algebraic data for $\bf{Ising}$}

There are three isomorphism classes of simple objects in $\bf{Ising}$, labeled $\{1,\sigma, \psi\}$. 

\begin{tabular}{ll}
Quantum dimensions & $d_1=1,\; d_{\sigma}=\sqrt{2},\;\hspace{14pt} d_{\psi}=1$ \\
Twists & $\theta_1=1,\; \theta_{\sigma}=e^{2\pi i/16},\; \theta_{\psi}=-1$ \\

$S$-matrix & $\frac{1}{2}\begin{pmatrix} 1 & \sqrt{2} & 1 \\ \sqrt{2} & 0 & -\sqrt{2} \\ 1 & -\sqrt{2} & 1 \\
\end{pmatrix}$ \\
$F$-symbols & $F^{\sigma\sigma\sigma}_{\sigma}=\frac{1}{\sqrt{2}} \begin{pmatrix}1 & 1 \\ 1 & -1 \end{pmatrix}$  \\
\end{tabular}
\subsection{Algebraic data for $\bf{Ising}\boxtimes \bf{Ising}$}
We also use the quantum dimensions and $S$-matrix of the doubled theory $\bf{Ising} \boxtimes \bf{Ising}$. We write objects in the bilayer category as $ab:= a\boxtimes b$ for simplicity. Under the Deligne product quantum dimensions simply multiply
$$d_{ab}=d_ad_b$$
and
the $S$-matrix respects the tensor product
$$S_{\bf{Ising} \boxtimes \bf{Ising}}=S_{\bf{Ising}}\otimes S_{\bf{Ising}}$$ but we nevertheless provide the values for convenience. The $S$-matrix is written with respect to the basis $\{11, 1\sigma, 1\psi, \sigma 1,\sigma \sigma, \sigma\psi, \psi 1, \psi \sigma, \psi \psi\}$. 

\begin{center}
\begin{tabular}{ll}
Quantum dimensions & $d_{11}=d_{1\psi}=d_{\psi 1}=d_{\psi\psi}=1, \;d_{\sigma 1}=d_{1\sigma}=d_{\sigma \psi}=d_{\psi \sigma}=\sqrt{2}, \; d_{\sigma\sigma}=2$ \\

$S$-matrix&  $\frac{1}{4}\begin{pmatrix} 1 & \sqrt{2} & 1 & \sqrt{2} & 2 & \sqrt{2} & 1 & \sqrt{2} & 1 \\ \sqrt{2} & 0 &-\sqrt{2} & 2 & 0 & -2 & \sqrt{2} & 0 & - \sqrt{2} \\
1 & -\sqrt{2} & 1 & \sqrt{2} & -2 &\sqrt{2} & 1 &-\sqrt{2} & 1 \\ \sqrt{2} &  2  & \sqrt{2} & 0 & 0 & 0 & -\sqrt{2} & -2 &-\sqrt{2} \\
2 & 0 & -2 & 0 & 0 & 0 & -2 & 0 & 2 \\ \sqrt{2} & -2 & \sqrt{2} & 0 & 0 & 0 & -\sqrt{2} & 2 & -\sqrt{2} \\ 1 & \sqrt{2} & 1 & -\sqrt{2} & -2 & -\sqrt{2} & 1 & \sqrt{2} & 1 \\ \sqrt{2} & 0 & -\sqrt{2} & -2 & 0 & 2 & \sqrt{2} & 0 & -\sqrt{2} \\ 1 & -\sqrt{2} & 1 & -\sqrt{2} & 2 &-\sqrt{2} & 1 & -\sqrt{2} & 1 \end{pmatrix}$ \\
\end{tabular}
\end{center}
\subsection{Algebraic data for $\left( \bf{Ising}\boxtimes \bf{Ising}\right)^{\times}_{\mathds{Z}_2}$}
While the full set of $\ZZ_2$-crossed $R$-symbols, $F$-, $U$-, and $\eta$- symbols is not yet known, we will only use the quantum dimensions and fusion rules (which are known) along with Conjecture 4.1 which we restate below for convenience.\\
\begin{tabular}{ll}
Quantum dimensions of defects & $d_{X_1}=d_{X_{\psi}}=2,\;d_{X_{\sigma}}=2\sqrt{2} $\\
Fusion rules & $X_1 \otimes X_1 = 11 \oplus \sigma\sigma \oplus \psi \psi$ \\
(see Section 4) &\\
\end{tabular}

Hereafter we write $X:=X_1$. 

\begin{conj} There exists a set of solutions to the $\mathds{Z}_2$-crossed consistency equations associated to  $\textbf{Ising}\boxtimes \textbf{Ising}$ such that
$$R^{XX}_{aa}=\theta_{a} \hspace{5pt} \text{ for } a \in \textbf{Ising}$$ and $$[F^{XXX}_X]=S_{\textbf{Ising}}= \frac{1}{2}\begin{pmatrix} 1 & \sqrt{2} & 1 \\ \sqrt{2} & 0 & -\sqrt{2} \\ 1 & -\sqrt{2} & 1 \\ \end{pmatrix}.$$
\end{conj}

\section{Graphical calculus}

We recall several moves involving the $G$-crossed  $R$-, $U$-, and $\eta$-symbols from \cite{BBCW} in the case where the fusion rules of $\left(\CC \right)^{\times}_G$ are multiplicity-free.
\begin{eqnarray}
R^{a_{ g} b_{ h}} \hspace{5pt}=\hspace{5pt}
\pspicture[shift=-0.75](-0.1,-0.4)(1.3,1.4)
\small
  \psset{linewidth=0.5pt,linecolor=black,arrowscale=1.5,arrowinset=0.15}
  \psline(0.96,0.05)(0.2,1)
  \psline{->}(0.96,0.05)(0.28,0.9)
  \psline(0.24,0.05)(1,1)
  \psline[border=2pt]{->}(0.24,0.05)(0.92,0.9)
  \rput[bl]{0}(-0.1,1.1){$a_{ g}$}
  \rput[br]{0}(1.2,1.1){$b_{ h}$}
  \rput[bl]{0}(-0.15,-0.05){$b_{ h}$}
  \rput[br]{0}(1.55,-0.1){$^{ \bar h} a_{ g}$}
  \endpspicture \hspace{5pt}=\hspace{5pt} \sum\limits_{c}\sqrt{\frac{d_{c_{gh}}}{d_{a_g}d_{b_h}}}
R_{c_{ gh}}^{a_{ g} b_{ h}}
 \pspicture[shift=-1](-0.3,-0.85)(1.8,1.3)
  \small
  \psset{linewidth=0.5pt,linecolor=black,arrowscale=1.5,arrowinset=0.15}
  \psline{->}(0.7,0)(0.7,0.45)
  \psline(0.7,0)(0.7,0.55)
  \psline(0.7,0.55) (0.25,1)
  \psline{->}(0.7,0.55)(0.3,0.95)
  \psline(0.7,0.55) (1.15,1)
  \psline{->}(0.7,0.55)(1.1,0.95)
  \rput[bl]{0}(0.1,0.2){$c_{ gh}$}
  \rput[br]{0}(1.4,1.05){$b_{ h}$}
  \rput[bl]{0}(-0.1,1.05){$a_{ g}$}
  \psline(0.7,0) (0.25,-0.45)
  \psline{-<}(0.7,0)(0.35,-0.35)
  \psline(0.7,0) (1.15,-0.45)
  \psline{-<}(0.7,0)(1.05,-0.35)
  \rput[br]{0}(1.7,-0.5){$^{ \bar h}a_{ g}$}
  \rput[bl]{0}(-0.1,-0.45){$b_{ h}$}
\scriptsize
  \endpspicture
  \end{eqnarray}
\begin{eqnarray}
\left(R^{ a_{\bf g} b_{\bf h}}\right)^{-1} =
\pspicture[shift=-0.75](-0.1,-0.4)(1.5,1.25)
\small
  \psset{linewidth=0.5pt,linecolor=black,arrowscale=1.5,arrowinset=0.15}
  \psline{->}(0.24,0.05)(0.92,0.9)
  \psline(0.24,0.05)(1,1)
  \psline(0.96,0.05)(0.2,1)
  \psline[border=2pt]{->}(0.96,0.05)(0.28,0.9)
  \rput[bl]{0}(0.0,1.1){$b_{h}$}
  \rput[br]{0}(1.3,1.05){$^{ \bar h}a_{ g}$}
  \rput[bl]{0}(0,-0.3){$a_{g}$}
  \rput[br]{0}(1.3,-0.25){$b_{ h}$}
  \endpspicture
\hspace{5pt}=\hspace{5pt}\sum\limits_{c}\sqrt{\frac{d_{c_{gh}}}{d_{a_g}d_{b_h}}}\left(
R_{c_{ gh}}^{a_{ g} b_{ h}}\right)^{-1}
 \pspicture[shift=-1](-0.1,-0.85)(1.5,1.4)
  \small
  \psset{linewidth=0.5pt,linecolor=black,arrowscale=1.5,arrowinset=0.15}
  \psline{->}(0.7,0)(0.7,0.45)
  \psline(0.7,0)(0.7,0.55)
  \psline(0.7,0.55) (0.25,1)
  \psline{->}(0.7,0.55)(0.3,0.95)
  \psline(0.7,0.55) (1.15,1)
  \psline{->}(0.7,0.55)(1.1,0.95)
  \rput[bl]{0}(0.1,0.2){$c_{gh}$}
  \rput[br]{0}(1.4,1.05){$^{ \bar h}a_{ g}$}
  \rput[bl]{0}(0,1.1){$b_{ h}$}
  \psline(0.7,0) (0.25,-0.45)
  \psline{-<}(0.7,0)(0.35,-0.35)
  \psline(0.7,0) (1.15,-0.45)
  \psline{-<}(0.7,0)(1.05,-0.35)
  \rput[br]{0}(1.4,-0.75){$b_{h}$}
  \rput[bl]{0}(0,-0.8){$a_{ g}$}
\scriptsize
  \endpspicture
\end{eqnarray}

\begin{eqnarray}
\psscalebox{1}{
\pspicture[shift=-1.7](-0.8,-0.8)(1.8,2.4)
  \small
  \psset{linewidth=0.5pt,linecolor=black,arrowscale=1.5,arrowinset=0.15}
  \psline{->}(0.7,0)(0.7,0.45)
  \psline(0.7,0)(0.7,0.55)
  \psline(0.7,0.55)(0.25,1)
  \psline(0.7,0.55)(1.15,1)	
  \psline(0.25,1)(0.25,2)
  \psline{->}(0.25,1)(0.25,1.9)
  \psline(1.15,1)(1.15,2)
  \psline{->}(1.15,1)(1.15,1.9)
  \psline[border=2pt](-0.65,0)(2.05,2)
  \psline{->}(-0.65,0)(0.025,0.5)
  \rput[bl]{0}(-0.4,0.6){$x_{ k}$}
  \rput[br]{0}(1.5,0.65){$^{\bar k}b$}
  \rput[bl]{0}(0.4,-0.4){$^{\bar k}c_{ gh}$}
  \rput[br]{0}(1.3,2.1){$b_{ h}$}
  \rput[bl]{0}(0.0,2.05){$a_{g}$}
 \scriptsize
  \endpspicture
}
=
U_{ k}\left(a, b; c\right)
\psscalebox{1}{
\pspicture[shift=-1.7](-0.8,-0.8)(1.8,2.4)
  \small
  \psset{linewidth=0.5pt,linecolor=black,arrowscale=1.5,arrowinset=0.15}
  \psline{->}(0.7,0)(0.7,0.45)
  \psline(0.7,0)(0.7,1.55)
  \psline(0.7,1.55)(0.25,2)
  \psline{->}(0.7,1.55)(0.3,1.95)
  \psline(0.7,1.55) (1.15,2)	
  \psline{->}(0.7,1.55)(1.1,1.95)
  \psline[border=2pt](-0.65,0)(2.05,2)
  \psline{->}(-0.65,0)(0.025,0.5)
  \rput[bl]{0}(-0.4,0.6){$x_{ k}$}
  \rput[bl]{0}(0.4,-0.4){$^{ \bar k}c_{gh}$}
  \rput[bl]{0}(0.15,1.2){$c_{ gh}$}
  \rput[br]{0}(1.3,2.1){$b_{h}$}
  \rput[bl]{0}(0.0,2.05){$a_{ g}$}
 \scriptsize
  \endpspicture
  }
  \end{eqnarray}
  \begin{eqnarray}
\psscalebox{1}{
\pspicture[shift=-1.7](-0.8,-0.8)(1.8,2.4)
  \small
  \psset{linewidth=0.5pt,linecolor=black,arrowscale=1.5,arrowinset=0.15}
  \psline{->}(0.7,0)(0.7,0.45)
  \psline(0.7,0)(0.7,1.55)
  \psline(0.7,1.55)(0.25,2)
  \psline{->}(0.7,1.55)(0.3,1.95)
  \psline(0.7,1.55) (1.15,2)	
  \psline{->}(0.7,1.55)(1.1,1.95)
  \psline[border=2pt](-0.65,2)(2.05,0)
  \psline{->}(2.05,0)(1.375,0.5)
  \rput[bl]{0}(1.4,0.6){$x_{ k}$}
  \rput[bl]{0}(0.4,-0.4){$^{ k}c_{ gh}$}
  \rput[bl]{0}(0.85,1.2){$c_{gh}$}
  \rput[br]{0}(1.3,2.1){$b_{h}$}
  \rput[bl]{0}(0.0,2.05){$a_{ g}$}
 \scriptsize
  \endpspicture
}
=
U_{ k}\left(\,^{ k}a,\,^{ k}b;\,^{k}c\right)
\psscalebox{1}{
\pspicture[shift=-1.7](-0.8,-0.8)(1.8,2.4)
  \small
  \psset{linewidth=0.5pt,linecolor=black,arrowscale=1.5,arrowinset=0.15}
  \psline{->}(0.7,0)(0.7,0.45)
  \psline(0.7,0)(0.7,0.55)
  \psline(0.7,0.55)(0.25,1)
  \psline(0.7,0.55)(1.15,1)	
  \psline(0.25,1)(0.25,2)
  \psline{->}(0.25,1)(0.25,1.9)
  \psline(1.15,1)(1.15,2)
  \psline{->}(1.15,1)(1.15,1.9)
  \psline[border=2pt](-0.65,2)(2.05,0)
  \psline{->}(2.05,0)(1.375,0.5)
  \rput[bl]{0}(1.4,0.6){$x_{ k}$}
  \rput[bl]{0}(0.4,-0.4){$^{k}c_{ gh}$}
  \rput[br]{0}(1.3,2.1){$b_{ h}$}
  \rput[bl]{0}(0.0,2.05){$a_{g}$}
  \rput[bl]{0}(-0.2,0.8){$^{ k}a$}
  \rput[bl]{0}(0.8,0.3){$^{ k}b$}
 \scriptsize
  \endpspicture
}
\qquad
\end{eqnarray}

\begin{eqnarray}
\psscalebox{1}{
\pspicture[shift=-1.7](-0.8,-0.8)(1.8,2.4)
  \small
  \psset{linewidth=0.5pt,linecolor=black,arrowscale=1.5,arrowinset=0.15}
  \psline(-0.65,0)(2.05,2)
  \psline[border=2pt](0.7,0.55)(0.25,1)
  \psline[border=2pt](1.15,1)(1.15,2)
  \psline(0.7,0.55)(1.15,1)	
  \psline{->}(0.7,0)(0.7,0.45)
  \psline(0.7,0)(0.7,0.55)
  \psline(0.25,1)(0.25,2)
  \psline{->}(0.25,1)(0.25,1.9)
  \psline{->}(1.15,1)(1.15,1.9)
  \psline{->}(-0.65,0)(0.025,0.5)
  \rput[bl]{0}(-0.5,0.6){$x_{ k}$}
  \rput[bl]{0}(0.5,1.3){$^{\bar g}x$}
  \rput[bl]{0}(1.6,2.1){$^{ \overline{gh}}x_{ k}$}
  \rput[bl]{0}(0.4,-0.35){$c_{ gh}$}
  \rput[br]{0}(1.3,2.1){$b_{ h}$}
  \rput[bl]{0}(0.0,2.05){$a_{ g}$}
 \scriptsize
  \endpspicture
}
=
\eta_{x}\left({ g},{ h}\right)
\psscalebox{1}{
\pspicture[shift=-1.7](-0.8,-0.8)(1.8,2.4)
  \small
  \psset{linewidth=0.5pt,linecolor=black,arrowscale=1.5,arrowinset=0.15}
  \psline(-0.65,0)(2.05,2)
  \psline[border=2pt](0.7,0)(0.7,1.55)
  \psline{->}(0.7,0)(0.7,0.45)
  \psline(0.7,1.55)(0.25,2)
  \psline{->}(0.7,1.55)(0.3,1.95)
  \psline(0.7,1.55) (1.15,2)	
  \psline{->}(0.7,1.55)(1.1,1.95)
  \psline{->}(-0.65,0)(0.025,0.5)
  \rput[bl]{0}(-0.5,0.6){$x_{ k}$}
  \rput[bl]{0}(1.6,2.1){$^{\overline{gh}}x_{ k}$}
  \rput[bl]{0}(0.4,-0.35){$c_{ gh}$}
  \rput[br]{0}(1.3,2.1){$b_{ h}$}
  \rput[bl]{0}(0.0,2.05){$a_{ g}$}
 \scriptsize
  \endpspicture
}
\end{eqnarray}
\begin{eqnarray}
\psscalebox{1}{
\pspicture[shift=-1.7](-0.8,-0.8)(1.8,2.4)
  \small
  \psset{linewidth=0.5pt,linecolor=black,arrowscale=1.5,arrowinset=0.15}
  \psline(2.05,0)(-0.65,2)
  \psline[border=2pt](0.7,0)(0.7,1.55)
  \psline{->}(0.7,0)(0.7,0.45)
  \psline(0.7,1.55)(0.25,2)
  \psline{->}(0.7,1.55)(0.3,1.95)
  \psline(0.7,1.55) (1.15,2)	
  \psline{->}(0.7,1.55)(1.1,1.95)
  \psline{->}(2.05,0)(1.375,0.5)
  \rput[bl]{0}(-0.8,2.1){$x_{ k}$}
  \rput[bl]{0}(1.4,0.6){$^{\overline{gh}} x_{ k}$}
  \rput[bl]{0}(0.4,-0.35){$c_{ gh}$}
  \rput[br]{0}(1.3,2.1){$b_{ h}$}
  \rput[bl]{0}(0.0,2.05){$a_{ g}$}
 \scriptsize
  \endpspicture
}
=
\eta_{x}\left({ g},{ h}\right)
\psscalebox{1}{
\pspicture[shift=-1.7](-0.8,-0.8)(1.8,2.4)
  \small
  \psset{linewidth=0.5pt,linecolor=black,arrowscale=1.5,arrowinset=0.15}
  \psline(2.05,0)(-0.65,2)
  \psline[border=2pt](0.7,0.55)(1.15,1)
  \psline[border=2pt](0.25,1)(0.25,2)
  \psline[border=2pt](1.15,1)(1.15,2)
  \psline(0.7,0.55)(0.25,1)
  \psline{->}(0.7,0)(0.7,0.45)
  \psline(0.7,0)(0.7,0.55)
  \psline(0.25,1)(0.25,2)
  \psline{->}(0.25,1)(0.25,1.9)
  \psline{->}(1.15,1)(1.15,1.9)
  \psline{->}(2.05,0)(1.375,0.5)
  \rput[bl]{0}(-0.8,2.1){$x_{ k}$}
  \rput[bl]{0}(1.4,0.6){$^{\overline{gh}} x_{ k}$}
  \rput[bl]{0}(0.5,1.3){$^{\bar g}x$}
  \rput[bl]{0}(0.4,-0.35){$c_{ gh}$}
  \rput[br]{0}(1.3,2.1){$b_{ h}$}
  \rput[bl]{0}(0.0,2.05){$a_{g}$}
 \scriptsize
  \endpspicture
}
\end{eqnarray}

\begin{equation}
\pspicture[shift=-1.0](-0.5,-0.9)(1.5,1.5)
\small
  \psset{linewidth=0.5pt,linecolor=black,arrowscale=1.5,arrowinset=0.15}
  \psline(0.75,-0.7)(0.75,-0.15)
  \psline(0.75,0.15)(0.75,1.3)
  \psellipse[linewidth=0.5pt,linecolor=black,border=0](0.4,0.3)(0.8,0.35)
  \psline{-<}(0.2,-0.027)(0.3,-0.04)
\psline[linewidth=0.5pt,linecolor=black,border=2.2pt,arrows=->,arrowscale=1.5,
arrowinset=0.15](0.75,0.1)(0.75,1.05)
\rput[bl]{0}(0.15,-0.4){$a_g$}
  \rput[tl]{0}(0.9,1.1){$b_h$}
\endpspicture
=\frac{S_{a_gb_h}}{S_{1b_h}}
\pspicture[shift=-1.0](0.4,-0.9)(1.5,1.5)
\small
  \psset{linewidth=0.5pt,linecolor=black,arrowscale=1.5,arrowinset=0.15}
  \psline(0.75,-0.7)(0.75,1.3)
\psline[linewidth=0.5pt,linecolor=black,border=2.2pt,arrows=->,arrowscale=1.5,
arrowinset=0.15](0.75,0.1)(0.75,1.05)
  \rput[tl]{0}(0.9,1.1){$b_h$}
 \endpspicture
 \end{equation}
\begin{eqnarray}
\pspicture[shift=-1.0](-0.5,-0.9)(1.5,1.5)
\small
  \psset{linewidth=0.5pt,linecolor=black,arrowscale=1.5,arrowinset=0.15}
  \psline(0.05,-0.7)(0.05,-0.15)
  \psline(0.05,0.15)(0.05,1.3)
  \psellipse[linewidth=0.5pt,linecolor=black,border=0](0.4,0.3)(0.8,0.35)
  \psline{-<}(0.2,-0.027)(0.3,-0.04)
\psline[linewidth=0.5pt,linecolor=black,border=2.2pt,arrows=->,arrowscale=1.5,arrowinset=0.15](0.05,0.1)(0.05,1.05)
\rput[bl]{0}(0.3,-0.4){$a_g$}
  \rput[tl]{0}(0.2,1.1){$b_h$}
\endpspicture
&=& \theta_{a_g}
\pspicture[shift=-1.0](-0.5,-0.9)(1.7,1.5)
\small
  \psset{linewidth=0.5pt,linecolor=black,arrowscale=1.5,arrowinset=0.15}
  \psarc(0.1, 0.2){0.4}{90}{270}
  \psarc(1.4,0.2){0.3}{-90}{90}
  \psbezier(0.1, -0.2)(0.4,-0.2)(1.0,0.5)(1.4,0.5)
\psline[border=2.1pt](1.1, -0.7)(1.1,0.2)
  \psbezier[linewidth=0.5pt, border=2.0pt](0.1, 0.6)(0.4,0.6)(1.0,-0.1)(1.4,-0.1)
  \psline[linewidth=.5pt,arrows=->, arrowscale=1.5, arrowinset=0.1](0.5,0.41)(0.4,0.48)
  \psline[linewidth=0.5pt,arrows=->, arrowscale=1.5, arrowinset=0.1, border=2.1pt](1.1, 0.2)(1.1,1.0)
  \psline(1.1,0.9)(1.1, 1.3)
\rput[bl]{0}(0.3,0.6){$a_g$}
  \rput[tl]{0}(1.2,1.1){$b_h$}
\endpspicture
= \theta_{a_g}
\pspicture[shift=-1.0](-0.5,-0.9)(1.8,1.8)
\small
  \psset{linewidth=0.5pt,linecolor=black,arrowscale=1.5,arrowinset=0.15}
  \psarc(0.1, 0.2){0.4}{90}{270}
  \psarc(1.4,0.2){0.3}{-90}{90}
  \psbezier(0.1, -0.2)(0.4,-0.2)(1.0,0.5)(1.4,0.5)
  \psbezier[linewidth=.5pt, border=2.0pt](0.1, 0.6)(0.4,0.6)(1.0,-0.1)(1.4,-0.1)
  \psline[linewidth=0.5pt,arrows=->, arrowscale=1.5, arrowinset=0.1](1.0,0.06)(0.9,0.12)
  \psline[border=2.1pt](0.4, -0.7)(0.4,0.35)
  \psline[linewidth=0.5pt,arrows=->, arrowscale=1.5, arrowinset=0.1](0.4, 0.55)(0.4,1.0)
  \psline(0.4,0.9)(0.4, 1.3)
\rput[bl]{0}(0.8,-0.3){$a_g$}
  \rput[tl]{0}(0.5,1.1){$b_h$}
\endpspicture
=
\pspicture[shift=-1.0](-0.5,-0.9)(1.3,1.5)
\small
  \psset{linewidth=0.5pt,linecolor=black,arrowscale=1.5,arrowinset=0.15}
  \psellipse[linewidth=0.5pt,linecolor=black,border=0](0.4,0.3)(0.8,0.35)
  \psline[linewidth=0.5pt,linecolor=black,border=2.2pt](0.75,-0.7)(0.75,0.5)
  \psline[linewidth=.5pt, arrows=->,arrowscale=1.5, arrowinset=0.15] (0.75,0.65) (0.75, 1.05)
  \psline{-<}(0.2,0.627)(0.3,0.64)
\psline(0.75,0.7)(0.75,1.1)
\psline (0.75, 1.0) (0.75, 1.3)
\rput[bl]{0}(-0.2,0.7){$a_g$}
  \rput[tl]{0}(0.9,1.1){$b_h$}
  \endpspicture
\label{}
\end{eqnarray}

\begin{equation}
  \pspicture[shift=-0.95](-0.2,-0.35)(1.2,1.75)
  \psarc[linewidth=.5pt,linecolor=black,border=0pt] (0.8,0.7){0.4}{120}{240}
  \psarc[linewidth=.5pt,linecolor=black,arrows=<-,arrowscale=1.4,
    arrowinset=0.15] (0.8,0.7){0.4}{165}{240}
  \psarc[linewidth=.5pt,linecolor=black,border=0pt] (0.4,0.7){0.4}{-60}{60}
  \psarc[linewidth=.5pt,linecolor=black,arrows=->,arrowscale=1.4,
    arrowinset=0.15] (0.4,0.7){0.4}{-60}{15}
  \psset{linewidth=.5pt,linecolor=black,arrowscale=1.5,arrowinset=0.15}
  \psline(0.6,1.05)(0.6,1.55)
  \psline{->}(0.6,1.05)(0.6,1.45)
  \psline(0.6,-0.15)(0.6,0.35)
  \psline{->}(0.6,-0.15)(0.6,0.25)
  \rput[bl]{0}(0,0.55){$a_g$}
  \rput[bl]{0}(0.94,0.55){$b_h$}
  \rput[bl]{0}(0.09,1.25){$c_{gh}$}
  \rput[bl]{0}(0.07,-0.05){$c_{gh}'$}
 \scriptsize

  \endpspicture
=\delta _{c_{gh} c_{gh} ^{\prime }}\sqrt{\frac{d_{a_g}d_{b_h}}{d_{c_{gh}}}}
  \pspicture[shift=-0.95](0.15,-0.35)(0.8,1.75)
  \small
  \psset{linewidth=.5pt,linecolor=black,arrowscale=1.5,arrowinset=0.15}
  \psline(0.6,-0.15)(0.6,1.55)
  \psline{->}(0.6,-0.15)(0.6,0.85)
  \rput[bl]{0}(0.75,1.25){$c_{gh}$}
  \endpspicture
\end{equation}%

\section{Calculation of matrix entries of the $T$-gate protocol}

\subsection{Qudit and protocol encoding}

In terms of the graphical calculus, the qubit is given by the $\mathds{C}$-span of the (normalized) fusion tree basis $|x\rangle$, which has the dual basis $\langle x |$.

\begin{minipage}{.5\textwidth}

$|x\rangle=$ 
\begin{tikzpicture}[scale=.5,baseline=-60,line width=.5]
\draw (1,-3/2)--(8,-5);
\draw (2,-3/2)--(2,-2);
\draw (3,-3/2)--(3,-5/2);
\draw (4,-3/2)--(4,-3);
\draw (5,-3/2)--(5,-7/2);
\draw (6,-3/2)--(6,-8/2);
\draw (7,-3/2)--(7,-9/2);
\draw (8,-3/2)--(8, -6);
\draw (1,-1) node {$\small \sigma 1$};
\draw (2,-1) node {$\small \sigma 1$};
\draw (3,-1) node {$\small \sigma 1$};
\draw (4,-1) node {$\small \sigma 1$};
\draw (5,-1) node {$\small X_1$};
\draw (6,-1) node {$\small X_1$};
\draw (7,-1) node {$\small X_1$};
\draw (8,-1) node {$\small X_1$};
\draw (2.25, -2.75) node {\small $x$};
\draw (3.5, -3.25) node {\small $\sigma 1$};
\draw(4.5,-3.70) node {\small $11$};
\draw(5.5,-4.25) node {\small $X_1$};
\draw(6.5,-4.75) node {\small $1$};
\draw(7.5,-5.25) node {\small $X_1$};
\draw (8,-6.5) node {\small $11$};
\end{tikzpicture}

\end{minipage}%
\begin{minipage}{.5\textwidth}
$\langle x |=$\scalebox{1}[-1]{
\begin{tikzpicture}[scale=.5,baseline=-47.5,line width=.5]
\draw (1,-3/2)--(8,-5);
\draw (2,-3/2)--(2,-2);
\draw (3,-3/2)--(3,-5/2);
\draw (4,-3/2)--(4,-3);
\draw (5,-3/2)--(5,-7/2);
\draw (6,-3/2)--(6,-8/2);
\draw (7,-3/2)--(7,-9/2);
\draw (8,-3/2)--(8, -6);
\draw (1,-1) node {\scalebox{1}[-1]{$\small \sigma 1$}};
\draw (2,-1) node {\scalebox{1}[-1]{$\small \sigma 1$}};
\draw (3,-1) node {\scalebox{1}[-1]{$\small \sigma 1$}};
\draw (4,-1) node {\scalebox{1}[-1]{$\small \sigma 1$}};
\draw (5,-1) node {\scalebox{1}[-1]{$\small X_1$}};
\draw (6,-1) node {\scalebox{1}[-1]{$\small X_1$}};
\draw (7,-1) node {\scalebox{1}[-1]{$\small X_1$}};
\draw (8,-1) node {\scalebox{1}[-1]{$\small X_1$}};
\draw (2.25, -2.75) node {\scalebox{1}[-1]{$\small x$}};
\draw (3.5, -3.25) node {\scalebox{1}[-1]{$\small \sigma 1$}};
\draw(4.5,-3.70) node {\scalebox{1}[-1]{$\small 11$}};
\draw(5.5,-4.25) node {\scalebox{1}[-1]{$\small X_1$}};
\draw(6.5,-4.75) node {\scalebox{1}[-1]{$\small 1$}};
\draw(7.5,-5.25) node {\scalebox{1}[-1]{$\small X_1$}};
\draw (8,-6.5) node {\scalebox{1}[-1]{$\small 11$}};
\end{tikzpicture}.}
\end{minipage}
In terms of diagrams the matrix entries $\langle y | T | x \rangle$ are given by stacking and taking a trace. Below we show the general diagram for the matrix entries alongside the shorthand notation that we use in what follows. In particular, we suppress the index on the symmetry defect, writing $X:=X_1$, and since all anyons live in a single layer, we write $a:=a1$ for $a$ a simple object in $\textbf{Ising}$. Moreover, since the total charge of the fusion trees are trivial, we suppress the labeling of the tracial strand.

As was mentioned in Section 4.3.2, the off diagonal matrix elements vanish by a Schur's Lemma argument, so it remains to calculate $\langle 1 | T | 1 \rangle$ and $\langle \psi | T | \psi \rangle$.

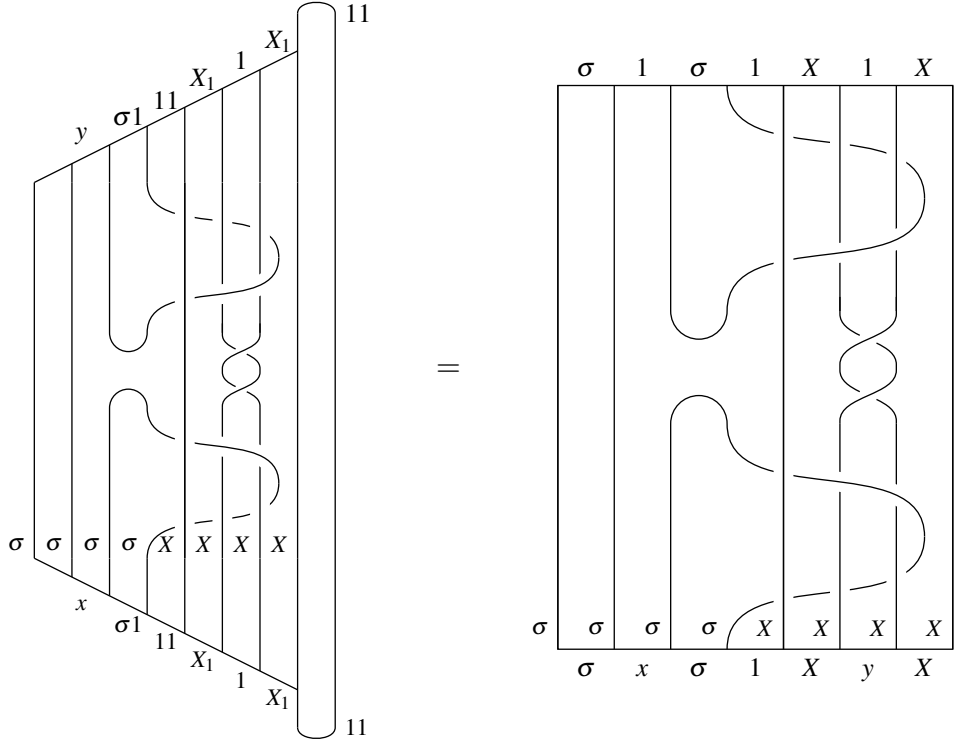
\begin{figure}[H]
\begin{minipage}{.45\textwidth}
\begin{tikzpicture}[scale=.5,baseline=-60,line width=.5]
\begin{scope}[shift={(0,14/2)}]
\draw (1,3/2)--(8,5);
\draw (2,3/2)--(2,2);
\draw (3,3/2)--(3,5/2);
\draw (4,3/2)--(4,3);
\draw (5,3/2)--(5,7/2);
\draw (6,3/2)--(6,8/2);
\draw (7,3/2)--(7,9/2);
\draw (8,3/2)--(8, 6);
\draw (2.25, 2.75) node {$\small y$};
\draw (3.5, 3.25) node {$\small \sigma 1$};
\draw(4.5,3.70) node {$\small 11$};
\draw(5.5,4.25) node {$\small X_1$};
\draw(6.5,4.75) node {$\small 1$};
\draw(7.5,5.25) node {$\small X_1$};
\draw (9,6) node[right] {$\small 11$};
\end{scope}

\draw(8,-6) to [out=-90, in=-90] (9,-6)--(9,13) to[out=90, in=90] (8,13);

\draw (1,-3/2)--(8,-5);
\draw (2,-3/2)--(2,-2);
\draw (3,-3/2)--(3,-5/2);
\draw (4,-3/2)--(4,-3);
\draw (5,-3/2)--(5,-7/2);
\draw (6,-3/2)--(6,-8/2);
\draw (7,-3/2)--(7,-9/2);
\draw (8,-3/2)--(8, -6);
\draw (2.25, -2.75) node {\small $x$};
\draw (3.5, -3.25) node {\small $\sigma 1$};
\draw(4.5,-3.70) node {\small $11$};
\draw(5.5,-4.25) node {\small $X_1$};
\draw(6.5,-4.75) node {\small $1$};
\draw(7.5,-5.25) node {\small $X_1$};
\draw (9,-6) node[right] {\small $11$};

\begin{scope}[shift={(0,-3/2)}]

\begin{scope}[shift={(5,6.25)}]
\braid a_1^{-1} a_1^{-1};
\end{scope}

\draw (4,0) to [out=90, in=-90] (7.5,2);
\draw (4,10) to [out=-90, in=90] (7.5,8);

\foreach \y in {6,7}
{
\draw[white,line width=7.5] (\y,0)--(\y,4);
\draw[white,line width=7.5] (\y,6.5)--(\y,10);
\draw (\y,0)--(\y,4);
\draw (\y,6)--(\y,10);
}

\draw[white, line width=7.5] (7.5,2) to [out=90, in=-90] (4,4) arc (0:180:.5)--(3,0);
\draw (7.5,2) to [out=90, in=-90] (4,4) arc (0:180:.5)--(3,0);

\draw[white, line width=7.5] (7.5,8) to [out=-90, in=90] (4,6) arc (0:-180:.5)--(3,10);
\draw (7.5,8) to [out=-90, in=90] (4,6) arc (0:-180:.5)--(3,10);


\draw[white,line width=7.5] (5,0)--(5,10);
\draw (5,0)--(5,10);
\foreach \x in {1,2,5,8}
{\draw (\x,0)--(\x,10);}

\foreach \x in {1,2,3,4}
{\draw (\x,.375) node [left] {\small $\sigma$};}
\foreach \y in {5,6,7,8}
{\draw (\y,.375) node [left] {\small $X$};}

\end{scope}
\end{tikzpicture}
\end{minipage}%
\begin{minipage}{.1\textwidth}
\btikz
\draw (5,0) node[left] {\Large $=$};
\etikz
\end{minipage}%
\begin{minipage}{.45\textwidth}
\btikz[line width=.5,scale=.75,baseline=70]

\begin{scope}[shift={(5,6.25)}]
\braid a_1^{-1} a_1^{-1};
\end{scope}

\draw (4,0) to [out=90, in=-90] (7.5,2);
\draw (4,10) to [out=-90, in=90] (7.5,8);

\foreach \y in {6,7}
{
\draw[white,line width=7.5] (\y,0)--(\y,4);
\draw[white,line width=7.5] (\y,6.5)--(\y,10);
\draw (\y,0)--(\y,4);
\draw (\y,6)--(\y,10);
}

\draw[white, line width=7.5] (7.5,2) to [out=90, in=-90] (4,4) arc (0:180:.5)--(3,0);
\draw (7.5,2) to [out=90, in=-90] (4,4) arc (0:180:.5)--(3,0);

\draw[white, line width=7.5] (7.5,8) to [out=-90, in=90] (4,6) arc (0:-180:.5)--(3,10);
\draw (7.5,8) to [out=-90, in=90] (4,6) arc (0:-180:.5)--(3,10);


\draw[white,line width=7.5] (5,0)--(5,10);
\draw (5,0)--(5,10);
\foreach \x in {1,2,5,8}
{\draw (\x,0)--(\x,10);}

\foreach \x in {1,2,3,4}
{\draw (\x,.375) node [left] {\small $\sigma$};}
\foreach \y in {5,6,7,8}
{\draw (\y,.375) node [left] {\small $X$};}

\draw (1,0) rectangle (8,10);
\draw (1.5,0) node[below] {$\phantom{X}\sigma\phantom{X}$};
\draw (2.5,0) node[below] {$\phantom{X}x\phantom{X}$};
\draw (3.5,0) node[below] {$\phantom{X}\sigma\phantom{X}$};
\draw (4.5,0) node[below] {$1$};
\draw (5.5,0) node[below] {$X$};
\draw (6.5,0) node[below] {$\phantom{X}y\phantom{X}$};
\draw (7.5,0) node[below] {$X$};
\draw (1.5,10) node[above] {$\phantom{X}\sigma\phantom{X}$};
\draw (2.5,10) node[above] {$\phantom{X}1\phantom{X}$};
\draw (3.5,10) node[above] {$\phantom{X}\sigma\phantom{X}$};
\draw (4.5,10) node[above] {$1$};
\draw (5.5,10) node[above] {$X$};
\draw (6.5,10) node[above] {$\phantom{X}1\phantom{X}$};
\draw (7.5,10) node[above] {$X$};
\etikz
\end{minipage}
\caption{On the left is the diagram representing the matrix entries of the protocol. On the right is the shorthand notation we use in the following section.} 
\end{figure}
\subsection{Calculation of $\langle 1 \mid T \mid 1 \rangle$}
We break down the calculation of the first matrix entry (up to normalization) into five steps. Below is a formula that shows the contributions from each step, before simplification.

$$\langle 1|T|1 \rangle=\underbrace{\frac{1}{d_{X}}\sum_c \sqrt{d_c} (R^{XX}_c)^2}_{\text{Step 1}} \underbrace{\Big ( \frac{S_{\sigma 1,c}}{S_{11,c}} \Big )^2}_{\text{Step 3}}\underbrace{\left | [F^{XXX}_{X}]_{c,11}\right |^2\frac{d_{X}^2}{d_c}}_{\text{Step 4}}\underbrace{d_{\sigma 1}d_{X}\sqrt{d_c}}_{\text {Step 5}}$$

The dashed rectangles in each equation indicate the region of the diagram where the graphical calculus is being applied.\\
\textbf{Step 1:} Unbraid the symmetry defects using equation (1) twice. 

\begin{eqnarray}
\btikz[line width=.5,scale=.5,baseline=70]

\begin{scope}[shift={(5,6.25)}]
\braid a_1^{-1} a_1^{-1};
\end{scope}

\draw (4,0) to [out=90, in=-90] (7.5,2);
\draw (4,10) to [out=-90, in=90] (7.5,8);

\foreach \y in {6,7}
{
\draw[white,line width=7.5] (\y,0)--(\y,4);
\draw[white,line width=7.5] (\y,6.5)--(\y,10);
\draw (\y,0)--(\y,4);
\draw (\y,6)--(\y,10);
}

\draw[white, line width=7.5] (7.5,2) to [out=90, in=-90] (4,4) arc (0:180:.5)--(3,0);
\draw (7.5,2) to [out=90, in=-90] (4,4) arc (0:180:.5)--(3,0);

\draw[white, line width=7.5] (7.5,8) to [out=-90, in=90] (4,6) arc (0:-180:.5)--(3,10);
\draw (7.5,8) to [out=-90, in=90] (4,6) arc (0:-180:.5)--(3,10);


\draw[white,line width=7.5] (5,0)--(5,10);
\draw (5,0)--(5,10);
\foreach \x in {1,2,5,8}
{\draw (\x,0)--(\x,10);}

\foreach \x in {1,2,3,4}
{\draw (\x,.375) node [left] {\small $\sigma$};}
\foreach \y in {5,6,7,8}
{\draw (\y,.375) node [left] {\small $X$};}
\stdlabels{1}{1}{1}{1}

\draw[red,dashed] (5.5,4) rectangle (7.5,6);
\draw (1,0) rectangle (8,10);
\etikz &=\sum_c\sqrt{\frac{d_{c}}{d_{X}d_{X}}}\left (R^{XX}_c\right)^2 \btikz[line width=.5,scale=.5,baseline=70]

\draw (6,6)--(6.5,5.25)--(7,6);

\draw (6,4)--(6.5,4.75)--(7,4);
\draw (6.5,4.75)--(6.5,5.25);
\draw (4,0) to [out=90, in=-90] (7.5,2);
\draw (4,10) to [out=-90, in=90] (7.5,8);

\foreach \y in {6,7}
{
\draw[white,line width=7.5] (\y,0)--(\y,3.5);
\draw[white,line width=7.5] (\y,6.5)--(\y,10);
\draw (\y,0)--(\y,4);
\draw (\y,6)--(\y,10);
}

\draw[white, line width=7.5] (7.5,2) to [out=90, in=-90] (4,4) arc (0:180:.5)--(3,0);
\draw (7.5,2) to [out=90, in=-90] (4,4) arc (0:180:.5)--(3,0);

\draw[white, line width=7.5] (7.5,8) to [out=-90, in=90] (4,6) arc (0:-180:.5)--(3,10);
\draw (7.5,8) to [out=-90, in=90] (4,6) arc (0:-180:.5)--(3,10);


\draw[white,line width=7.5] (5,0)--(5,10);
\draw (5,0)--(5,10);
\foreach \x in {1,2,5,8}
{\draw (\x,0)--(\x,10);}

\foreach \x in {1,2,3,4}
{\draw (\x,.375) node [left] {\small $\sigma$};}
\foreach \y in {5,6,7,8}
{\draw (\y,.375) node [left] {\small $X$};}
\stdlabels{1}{1}{1}{1}
\draw[red,dashed] (5.5,4) rectangle (7.5,6);
\draw (6.5,5) node[right] {$c$};
\draw (1,0) rectangle (8,10);

\etikz
\end{eqnarray}

Note that in this case there is a loop labeled by $\sigma$ that we can replace with $d_{\sigma}$ right away. However, we have left it unsimplified so that the same picture applies for the calculation of $\langle \psi | T | \psi \rangle$.

\textbf{Step 2:} Slide the $\sigma$ loops under the first defect charge line so that they encircle pairs of defect charge lines using equations (3)-(6). The $U$- and $\eta$- symbols this introduces all cancel, so that the anyons and defects can be moved past each other at no cost.

\begin{eqnarray}
\btikz[line width=.5,scale=.5,baseline=70]

\draw (6,6)--(6.5,5.25)--(7,6);

\draw (6,4)--(6.5,4.75)--(7,4);
\draw (6.5,4.75)--(6.5,5.25);
\draw (4,0) to [out=90, in=-90] (7.5,2);
\draw (4,10) to [out=-90, in=90] (7.5,8);

\foreach \y in {6,7}
{
\draw[white,line width=7.5] (\y,0)--(\y,3.5);
\draw[white,line width=7.5] (\y,6.5)--(\y,10);
\draw (\y,0)--(\y,4);
\draw (\y,6)--(\y,10);
}

\draw[white, line width=7.5] (7.5,2) to [out=90, in=-90] (4,4) arc (0:180:.5)--(3,0);
\draw (7.5,2) to [out=90, in=-90] (4,4) arc (0:180:.5)--(3,0);

\draw[white, line width=7.5] (7.5,8) to [out=-90, in=90] (4,6) arc (0:-180:.5)--(3,10);
\draw (7.5,8) to [out=-90, in=90] (4,6) arc (0:-180:.5)--(3,10);


\draw[white,line width=7.5] (5,0)--(5,10);
\draw (5,0)--(5,10);
\foreach \x in {1,2,5,8}
{\draw (\x,0)--(\x,10);}

\foreach \x in {1,2,3,4}
{\draw (\x,.375) node [left] {\small $\sigma$};}
\foreach \y in {5,6,7,8}
{\draw (\y,.375) node [left] {\small $X$};}
\stdlabels{1}{1}{1}{1}
\draw (6.5,5) node[right] {$c$};
\draw[dashed,red] (2.5,5.25) rectangle (7.75,10);
\draw[dashed,red] (2.5,0) rectangle (7.75,4.75);
\draw (1,0) rectangle (8,10);

\etikz
&=
\btikz[line width=.5,scale=.5,baseline=70]

\draw (6,6)--(6.5,5.25)--(7,6);

\draw (6,4)--(6.5,4.75)--(7,4);
\draw (6.5,4.75)--(6.5,5.25);


\draw[white,line width=7.5] (5,0)--(5,10);
\draw (5,0)--(5,10);
\foreach \x in {1,2,5,8}
{\draw (\x,0)--(\x,10);}

\foreach \x in {1,2}
{\draw (\x,.375) node [left] {\small $\sigma$};}
\foreach \y in {5,6,7,8}
{\draw (\y,.375) node [left] {\small $X$};}

\begin{scope}[ >=stealth,decoration={markings, mark=at position .6 with {\arrow{>}}}]
\draw[white, line width=7.5] (6.5,3) circle [x radius=1, y radius=.5, start angle = 0, end angle = 180];
\draw[postaction={decorate}] (6.5,3.125) circle [x radius=1, y radius=.5];

\draw[white, line width=7.5] (6.5,6.875) circle [x radius=1, y radius=.5, start angle = 0, end angle = -180];
\draw[postaction={decorate}] (6.5,6.875) circle [x radius=1, y radius=.5];

\end{scope}
\foreach \y in {6,7}{
\draw[white,line width=7.5] (\y,0)--(\y,3.25);
\draw[white,line width=7.5] (\y,6.75)--(\y,10);
\draw (\y,0)--(\y,3.25);
\draw (\y, 3.75)--(\y,4);
\draw (\y,6.75)--(\y,10);
\draw (\y, 6.25)--(\y,6);}

\draw (1.5,0) node[below] {$\phantom{X}\sigma\phantom{X}$};
\draw (3.5,0) node[below] {$\phantom{X}1\phantom{X}$};
\draw (5.5,0) node[below] {$X$};
\draw (6.5,0) node[below] {$\phantom{X}1\phantom{X}$};
\draw (7.5,0) node[below] {$X$};
\draw (1.5,10) node[above] {$\phantom{X}\sigma\phantom{X}$};
\draw (3.5,10) node[above] {$\phantom{X}1\phantom{X}$};
\draw (5.5,10) node[above] {$X$};
\draw (6.5,10) node[above] {$\phantom{X}1\phantom{X}$};
\draw (7.5,10) node[above] {$X$};
\draw (5.675, 6.125) node {\small $\sigma$};
\draw (5.675, 2.125) node {\small $\sigma$};
\draw (6.5,5) node[right] {$c$};
\draw[thick] (1,0) rectangle (8,10);
\draw[red,dashed] (5.25,1.75) rectangle (7.75,4);
\draw[red,dashed] (5.25,7.75) rectangle (7.75,5.5);
\draw (1,0) rectangle (8,10);

\etikz
\end{eqnarray}

\textbf{Step 3:} Slide the loops over the trivalent vertices so that they encircle the fusion channel labeled by $c$, flip the tilt of the bottom loop using (6), then use the loop removal relation (7) twice. \\
In general with loops encircling charge lines in the $G$-crossed graphical calculus we need to pay careful attention to the relative positions of charge lines, see equations (306) and (307) in \cite{BBCW}. However, in this case the subdiagram involves only anyons, so we can use the graphical calculus of a UMTC and ignore them. Nevertheless, a careful account of the factors using the full $G$-extended calculus gives the same result: any factors of $U$- and $\eta$- cancel. 
\begin{eqnarray}
\btikz[line width=.5,scale=.5,baseline=70]

\draw (6,6)--(6.5,5.25)--(7,6);
\draw (6,4)--(6.5,4.75)--(7,4);
\draw (6.5,4.75)--(6.5,5.25);


\draw[white,line width=7.5] (5,0)--(5,10);
\draw (5,0)--(5,10);
\foreach \x in {1,2,5,8}
{\draw (\x,0)--(\x,10);}

\foreach \x in {1,2}
{\draw (\x,.375) node [left] {\small $\sigma$};}
\foreach \y in {5,6,7,8}
{\draw (\y,.375) node [left] {\small $X$};}
\labelstwo{1}{1}
\draw (6.5,5) node[right] {$c$};

\begin{scope}[ >=stealth,decoration={markings, mark=at position .6 with {\arrow{>}}}]
\draw[white, line width=7.5] (6.5,3) circle [x radius=1, y radius=.5, start angle = 0, end angle = 180];
\draw[postaction={decorate}] (6.5,3.125) circle [x radius=1, y radius=.5];

\draw[white, line width=7.5] (6.5,6.875) circle [x radius=1, y radius=.5, start angle = 0, end angle = -180];
\draw[postaction={decorate}] (6.5,6.875) circle [x radius=1, y radius=.5];
\end{scope}
\foreach \y in {6,7}{
\draw[white,line width=7.5] (\y,0)--(\y,3.25);
\draw[white,line width=7.5] (\y,6.75)--(\y,10);
\draw (\y,0)--(\y,3.25);
\draw (\y, 3.75)--(\y,4);
\draw (\y,6.75)--(\y,10);
\draw (\y, 6.25)--(\y,6);}
\draw (5.675, 6.125) node {\small $\sigma$};
\draw (5.675, 2.125) node {\small $\sigma$};

\draw[thick] (1,0) rectangle (8,10);
\draw[red,dashed] (5.25,1.75) rectangle (7.75,7.75);
\etikz &=& \btikz[line width=.5,scale=.5,baseline=70]

\draw (6,8)--(6.5,7.25)--(7,8);
\draw (6,2)--(6.5,2.75)--(7,2);
\draw (6.5,2.75)--(6.5,7.25);

\foreach \x in {6,7}{
\draw  (\x,2)--(\x,0);
\draw (\x,8)--(\x,10);
}

\foreach \x in {1,2,5,8}
{\draw (\x,0)--(\x,10);}

\foreach \x in {1,2}
{\draw (\x,.375) node [left] {\small $\sigma$};}
\foreach \y in {5,6,7,8}
{\draw (\y,.375) node [left] {\small $X$};}
\labelstwo{1}{1}

\begin{scope}[ >=stealth,decoration={markings, mark=at position .6 with {\arrow{>}}}]
\draw[white, line width=7.5] (6.5,5.875) circle [x radius=1, y radius=.5, start angle = 0, end angle = -180];
\draw[postaction={decorate}] (6.5,5.875) circle [x radius=1, y radius=.5];

\draw[white, line width=7.5] (6.5,6)--(6.5,7);
\draw (6.5,6)--(6.5,7);

\draw[white, line width=7.5] (6.5,4.125) circle [x radius=1, y radius=.5, start angle = 0, end angle = 180];
\draw[postaction={decorate}]  (6.5,4.125) circle [x radius=1, y radius=.5];
\draw[white, line width=7.5] (6.5,4)--(6.5,3);
\draw (6.5,3)--(6.5,4);
\end{scope}
\draw[thick] (1,0) rectangle (8,10);
\draw[red,dashed] (5.25,2.5) rectangle (7.75,7);
\draw (6.5,5) node[right] {$c$};
\draw (5.675, 5.125) node {\small $\sigma$};
\draw (5.675, 3.125) node {\small $\sigma$};
\etikz \\
 & = \normalsize{ \left( \frac{S_{\sigma1,c}}{S_{11,c}}\right )^2} &\btikz[line width=.5,scale=.5,baseline=70]

\draw (6,8)--(6.5,7.25)--(7,8);
\draw (6,2)--(6.5,2.75)--(7,2);
\draw (6.5,2.75)--(6.5,7.25);

\foreach \x in {6,7}{
\draw  (\x,2)--(\x,0);
\draw (\x,8)--(\x,10);
}

\foreach \x in {1,2,5,8}
{\draw (\x,0)--(\x,10);}

\foreach \x in {1,2}
{\draw (\x,.375) node [left] {\small $\sigma$};}
\foreach \y in {5,6,7,8}
{\draw (\y,.375) node [left] {\small $X$};}
\labelstwo{1}{1}
\draw (6.5,5) node[right] {$c$};
\draw(1,0) rectangle (8,10);
\draw[red,dashed] (5.25,3) rectangle (7.75,7);

\etikz
\end{eqnarray}

\textbf{Step 4:} Perform a sequence of $F$-moves to change the basis states, then use the bubble-popping relation (9).

\begin{eqnarray}
\btikz[line width=.5,scale=.45,baseline=65]

\draw (6,8)--(6.5,7.25)--(7,8);
\draw (6,2)--(6.5,2.75)--(7,2);
\draw (6.5,2.75)--(6.5,7.25);

\foreach \x in {6,7}{
\draw  (\x,2)--(\x,0);
\draw (\x,8)--(\x,10);
}

\foreach \x in {1,2,5,8}
{\draw (\x,0)--(\x,10);}

\foreach \x in {1,2}
{\draw (\x,.375) node [left] {\small $\sigma$};}
\foreach \y in {5,6,7,8}
{\draw (\y,.375) node [left] {\small $X$};}
\labelstwo{1}{1}
\draw (6.5,5) node[right] {$c$};
\draw(1,0) rectangle (8,10);
\draw[red,dashed] (5.25,-1) rectangle (7.75,11);

\etikz =\sum_{d,d'} \left [F^{XXX}_X \right ]_{d,11} \left [F^{XXX}_X \right ]^*_{d',11} 
&\btikz[line width=.5,scale=.45,baseline=65]

\draw (6,8)--(6.5,7.25)--(7,8);
\draw (6,8)--(6.5,8.75)--(7,8);

\draw (6,2)--(6.5,2.75)--(7,2);
\draw (6,2)--(6.5,1.25)--(7,2);

\draw (6.5,2.75)--(6.5,7.25);
\draw (6.5,1.25)--(6.5,0);
\draw (6.5,8.75)--(6.5,10);


\foreach \x in {1,2,5,8}
{\draw (\x,0)--(\x,10);}

\foreach \x in {1,2}
{\draw (\x,.375) node [left] {\small $\sigma$};}
\foreach \y in {5,8}
{\draw (\y,.375) node [left] {\small $X$};}
\draw (6.5,5) node[right] {$c$};
\draw (6.5,1) node[left] {$d$};
\draw (6.5,9) node[left] {$d'$};
\labelsthree
\draw(1,0) rectangle (8,10);
\draw[red,dashed] (5.25,-1) rectangle (7.75,11);

\etikz \\
 =\sum_{d,d'} \left [F^{XXX}_X \right ]_{d,11} \left [F^{XXX}_X \right ]^*_{d',11} \left ( \sqrt{\frac{d_Xd_X}{d_c}}\right)^2 \delta_{dc}\delta_{d'c} &\btikz[line width=.5,scale=.45,baseline=65]

\draw (6.5,10)--(6.5,0);


\foreach \x in {1,2,5,8}
{\draw (\x,0)--(\x,10);}

\foreach \x in {1,2}
{\draw (\x,.375) node [left] {\small $\sigma$};}
\foreach \y in {5,8}
{\draw (\y,.375) node [left] {\small $X$};}
\draw (6.5,5) node[right] {$c$};
\labelsthree
\draw(1,0) rectangle (8,10);
\draw[red,dashed] (5.25,-1) rectangle (7.75,11);

\etikz
\end{eqnarray}

\textbf{Step 5:} The diagram that remains contributes a scalar factor that can be calculated using equations (4) and (5) repeatedly until the empty diagram is obtained.

\begin{eqnarray}
\btikz[line width=.5,scale=.45,baseline=65]

\draw (6.5,10)--(6.5,0);


\foreach \x in {1,2,5,8}
{\draw (\x,0)--(\x,10);}

\foreach \x in {1,2}
{\draw (\x,.375) node [left] {\small $\sigma$};}
\foreach \y in {5,8}
{\draw (\y,.375) node [left] {\small $X$};}
\draw (6.5,5) node[right] {$c$};
\labelsthree
\draw(1,0) rectangle (8,10);

\etikz &=& d_{\sigma 1}d_{X}\sqrt{d_c}
\end{eqnarray}

Putting these steps together, we arrive at the formula

\begin{eqnarray} \frac{1}{d_{X}}\sum_c \sqrt{d_c} (R^{XX}_c)^2 \Big ( \frac{S_{\sigma 1,c}}{S_{11,c}} \Big )^2\left |[F^{XXX}_{X}]_{c,11}\right |^2\frac{d_{X}^2}{d_c}d_{\sigma 1}d_{X}\sqrt{d_c} \\
= d_{\sigma 1} d_{X}^2 \sum_c  (R^{XX}_c)^2 \left |[F^{XXX}_{X}]_{c,11}^2\right |^2\Big ( \frac{S_{\sigma 1,c}}{S_{11,c}} \Big )^2
\end{eqnarray}

\subsection{Calculation of $\langle \psi | T | \psi \rangle$}
The first step of resolving the defect braiding is  identical to that of the previous subsection for $\langle 1 | T | 1 \rangle$. However, when the fusion channel of the four monolayer Ising anyons $\sigma 1$ is given by $\psi 1$, the loops are not free to slide under the defect line as in Step 2. Instead, we perform an $F$-move in the middle of the diagram.

\begin{eqnarray}
 \btikz[line width=.5,scale=.5,baseline=70]

\draw (6,6)--(6.5,5.25)--(7,6);

\draw (6,4)--(6.5,4.75)--(7,4);
\draw (6.5,4.75)--(6.5,5.25);
\draw (4,0) to [out=90, in=-90] (7.5,2);
\draw (4,10) to [out=-90, in=90] (7.5,8);

\foreach \y in {6,7}
{
\draw[white,line width=7.5] (\y,0)--(\y,3.5);
\draw[white,line width=7.5] (\y,6.5)--(\y,10);
\draw (\y,0)--(\y,4);
\draw (\y,6)--(\y,10);
}

\draw[white, line width=7.5] (7.5,2) to [out=90, in=-90] (4,4) arc (0:180:.5)--(3,0);
\draw (7.5,2) to [out=90, in=-90] (4,4) arc (0:180:.5)--(3,0);

\draw[white, line width=7.5] (7.5,8) to [out=-90, in=90] (4,6) arc (0:-180:.5)--(3,10);
\draw (7.5,8) to [out=-90, in=90] (4,6) arc (0:-180:.5)--(3,10);


\draw[white,line width=7.5] (5,0)--(5,10);
\draw (5,0)--(5,10);
\foreach \x in {1,2,5,8}
{\draw (\x,0)--(\x,10);}

\foreach \x in {1,2,3,4}
{\draw (\x,.375) node [left] {\small $\sigma$};}
\foreach \y in {5,6,7,8}
{\draw (\y,.375) node [left] {\small $X$};}
\draw (6.5,5) node[right] {$c$};
\draw[dashed] (5.75,3.0)--(5.75,7);
\draw[red,dashed] (5.25,1.25) rectangle (7.75,8.75);
\draw (1,0) rectangle (8,10);
\stdlabels{\psi}{1}{\psi}{1}
\etikz
&=&\sum_d [F^{\sigma 1 \sigma 1 \sigma 1}_{\sigma 1}]_{d,11}
 \btikz[line width=.5,scale=.5,baseline=70]

\draw (6,8.5)--(6.5,7.75)--(7,8.5);
\draw (6,1.5)--(6.5,2.25)--(7,1.5);
\draw (6.5,2.25)--(6.5,7.75);
\draw[white, line width=7.5] (4,0) to [out=90, in=-90] (7.5,2) to [out=90,in=-60] (6,5);
\draw (4,0) to [out=90, in=-90] (7.5,2) to [out=90,in=-60] (6,5);
\draw[white, line width=7.5] (4,10) to [out=-90, in=90] (7.5,8) to [out=-90, in=60] (6,5);

\draw (4,10) to [out=-90, in=90] (7.5,8) to [out=-90, in=60] (6,5);
\draw[blue] (6,5)--(5.5,5);
\draw (5.75,5) node[below] {$d$};

\draw (5.5,5) to [out=120, in=-90] (3,10);
\draw (5.5,5) to [out=240, in=90] (3,0);
\foreach \y in {6,7}
{
\draw[white,line width=7.5] (\y,0)--(\y,1.5);
\draw[white,line width=7.5] (\y,8.5)--(\y,10);
\draw (\y,0)--(\y,1.5);
\draw (\y,8.5)--(\y,10);
}


\draw[white,line width=7.5] (5,0)--(5,10);
\draw (5,0)--(5,10);
\foreach \x in {1,2,5,8}
{\draw (\x,0)--(\x,10);}

\foreach \x in {1,2,3,4}
{\draw (\x,.375) node [left] {\small $\sigma$};}
\foreach \y in {5,6,7,8}
{\draw (\y,.375) node [left] {\small $X$};}
\draw (6.5,5) node[right] {$c$};
\draw[red,dashed] (5.25,1.25) rectangle (7.75,8.75);
\draw (1,0) rectangle (8,10);
\stdlabels{\psi}{1}{\psi}{1}
\etikz
\end{eqnarray}

When $d$ is the vacuum channel, $d=11$, we can perform a sequence of sliding moves to resolve all crossings by the same arguments as in the previous subsection. 

\begin{eqnarray}
\btikz[line width=.5,scale=.5,baseline=70]

\draw (6,8.5)--(6.5,7.75)--(7,8.5);
\draw (6,1.5)--(6.5,2.25)--(7,1.5);
\draw (6.5,2.25)--(6.5,7.75);
\draw[white, line width=7.5] (4,0) to [out=90, in=-90] (7.5,2) to [out=90,in=-90] (6,5);
\draw (4,0) to [out=90, in=-90] (7.5,2) to [out=90,in=-90] (6,5);
\draw[white, line width=7.5] (4,10) to [out=-90, in=90] (7.5,8) to [out=-90, in=90] (6,5);
\draw (4,10) to [out=-90, in=90] (7.5,8) to [out=-90, in=90] (6,5);

\draw (5.5,5) to [out=90, in=-90] (3,10);
\draw (5.5,5) to [out=-90, in=90] (3,0);
\foreach \y in {6,7}
{
\draw[white,line width=7.5] (\y,0)--(\y,1.5);
\draw[white,line width=7.5] (\y,8.5)--(\y,10);
\draw (\y,0)--(\y,1.5);
\draw (\y,8.5)--(\y,10);
}




\draw[white,line width=7.5] (5,0)--(5,10);
\draw (5,0)--(5,10);
\foreach \x in {1,2,5,8}
{\draw (\x,0)--(\x,10);}

\foreach \x in {1,2,3,4}
{\draw (\x,.375) node [left] {\small $\sigma$};}
\foreach \y in {5,6,7,8}
{\draw (\y,.375) node [left] {\small $X$};}
\draw (6.5,5) node[right] {$c$};
\stdlabels{\psi}{1}{\psi}{1}
\draw (1,0) rectangle (8,10);
\etikz
&=&
\btikz[line width=.5,scale=.5,baseline=70]

\draw (6,8)--(6.5,7.25)--(7,8);
\draw (6,2)--(6.5,2.75)--(7,2);
\draw (6.5,2.75)--(6.5,7.25);

\foreach \x in {6,7}{
\draw  (\x,2)--(\x,0);
\draw (\x,8)--(\x,10);
}

\foreach \x in {1,2,5,8}
{\draw (\x,0)--(\x,10);}

\foreach \x in {1,2}
{\draw (\x,.375) node [left] {\small $\sigma$};}
\foreach \y in {5,6,7,8}
{\draw (\y,.375) node [left] {\small $X$};}
\draw (6.5,5) node[right] {$c$};
\labelstwo{\psi}{\psi}
\draw(1,0) rectangle (8,10);
\draw[red,dashed] (5.25,-1) rectangle (7.75,11);
\draw (1,0) rectangle (8,10);
\etikz
\end{eqnarray}
Now the righthand side can be resolved according as with steps 4 and 5 for $\langle 1 | T | 1 \rangle$, with the minor difference that the fusion channel of the monolayer Ising anyons is given by $\psi 1$. 

However, when $d=\psi 1$ we must use a more involved sequence of moves to resolve the diagram.

\begin{eqnarray}
 \btikz[line width=.5,scale=.5,baseline=70]

\draw (6,8.5)--(6.5,7.75)--(7,8.5);
\draw (6,1.5)--(6.5,2.25)--(7,1.5);
\draw (6.5,2.25)--(6.5,7.75);
\draw[white, line width=7.5] (4,0) to [out=90, in=-90] (7.5,2) to [out=90,in=-60] (6,5);
\draw (4,0) to [out=90, in=-90] (7.5,2) to [out=90,in=-60] (6,5);
\draw[white, line width=7.5] (4,10) to [out=-90, in=90] (7.5,8) to [out=-90, in=60] (6,5);

\draw (4,10) to [out=-90, in=90] (7.5,8) to [out=-90, in=60] (6,5);
\draw[blue] (6,5)--(5.5,5);
\draw (5.75,5) node[below] {$\psi$};

\draw (5.5,5) to [out=120, in=-90] (3,10);
\draw (5.5,5) to [out=240, in=90] (3,0);
\foreach \y in {6,7}
{
\draw[white,line width=7.5] (\y,0)--(\y,1.5);
\draw[white,line width=7.5] (\y,8.5)--(\y,10);
\draw (\y,0)--(\y,1.5);
\draw (\y,8.5)--(\y,10);
}


\draw[white,line width=7.5] (5,0)--(5,10);
\draw (5,0)--(5,10);
\foreach \x in {1,2,5,8}
{\draw (\x,0)--(\x,10);}

\foreach \x in {1,2,3,4}
{\draw (\x,.375) node [left] {\small $\sigma$};}
\foreach \y in {5,6,7,8}
{\draw (\y,.375) node [left] {\small $X$};}
\draw (1,0) rectangle (8,10);
\draw[red,dashed] (5.25,3.5) rectangle (7.75,6.5);
\draw (6.5,5) node[right] {$c$};
\stdlabels{\psi}{1}{\psi}{1}
\etikz
&=& \btikz[line width=.5,scale=.5,baseline=70]

\draw (6,8.5)--(6.5,7.75)--(7,8.5);
\draw (6,1.5)--(6.5,2.25)--(7,1.5);
\draw (6.5,2.25)--(6.5,7.75);

\draw[white, line width=7.5] (4,0) to [out=90, in=-90] (7.5,2) to [out=90,in=-60] (7,5);
\draw (4,0) to [out=90, in=-90] (7.5,2) to [out=90,in=-60] (7,5);
\draw[white, line width=7.5] (4,10) to [out=-90, in=90] (7.5,8) to [out=-90, in=60] (7,5);

\draw (4,10) to [out=-90, in=90] (7.5,8) to [out=-90, in=60] (7,5);
\draw[white, line width=7.5] (7,5)--(5.5,5);
\draw[blue] (7,5)--(5.5,5);
\draw (6.25,5) node[below] {$\psi$};

\draw (5.5,5) to [out=120, in=-90] (3,10);
\draw (5.5,5) to [out=240, in=90] (3,0);
\foreach \y in {6,7}
{
\draw[white,line width=7.5] (\y,0)--(\y,1.5);
\draw[white,line width=7.5] (\y,8.5)--(\y,10);
\draw (\y,0)--(\y,1.5);
\draw (\y,8.5)--(\y,10);
}




\draw[white,line width=7.5] (5,0)--(5,10);
\draw (5,0)--(5,10);
\foreach \x in {1,2,5,8}
{\draw (\x,0)--(\x,10);}

\foreach \x in {1,2,3,4}
{\draw (\x,.375) node [left] {\small $\sigma$};}
\foreach \y in {5,6,7,8}
{\draw (\y,.375) node [left] {\small $X$};}
\draw (6.5,2.5) node[right] {$c$};
\stdlabels{\psi}{1}{\psi}{1}
\draw[red,dashed] (5.25,3.5) rectangle (7.75,6.5);
\draw (1,0) rectangle (8,10);
\etikz \\
= \btikz[line width=.5,scale=.5,baseline=70]

\draw (6,8.5)--(6.5,7.75)--(7,8.5);
\draw (6,1.5)--(6.5,2.25)--(7,1.5);
\draw (6.5,2.25)--(6.5,7.75);

\draw[white, line width=7.5,looseness=1] (4,0) to [out=90, in=-90] (7.5,2) ;
\draw[white, line width=7.5,looseness=2.5]  (7.5,2) to [out=90,in=60] (5.5,5.5);

\draw[looseness=1] (4,0) to [out=90, in=-90] (7.5,2);
\draw[looseness=2.5] (7.5,2) to [out=90,in=60] (5.5,5.5);
\draw[white, line width=7.5,looseness=.5] (4,10) to [out=-90, in=150] (5.5,5.5);


\draw[looseness=.5] (4,10) to [out=-90, in=150] (5.5,5.5);

\draw[white, line width=7.5] (7,5)--(5.5,5);
\draw[blue] (5.5,5) to [out=0,in=-90] (7.25,5.25) to [out=90, in=0] (5.5,5.5);

\draw[white,line width=7.5] (6.5,5)--(6.5,7);
\draw (6.5,5)--(6.5,7.5);

\draw (6.25,5) node[below] {$\psi$};

\draw (5.5,5) to [out=120, in=-90] (3,10);
\draw (5.5,5) to [out=240, in=90] (3,0);
\foreach \y in {6,7}
{
\draw[white,line width=7.5] (\y,0)--(\y,1.5);
\draw[white,line width=7.5] (\y,8.5)--(\y,10);
\draw (\y,0)--(\y,1.5);
\draw (\y,8.5)--(\y,10);
}




\draw[white,line width=7.5] (5,0)--(5,10);
\draw (5,0)--(5,10);
\foreach \x in {1,2,5,8}
{\draw (\x,0)--(\x,10);}

\foreach \x in {1,2,3,4}
{\draw (\x,.375) node [left] {\small $\sigma$};}
\foreach \y in {5,6,7,8}
{\draw (\y,.375) node [left] {\small $X$};}
\draw (6.5,2.5) node[right] {$c$};
\stdlabels{\psi}{1}{\psi}{1}
\draw[red,dashed] (5.25,3.5) rectangle (7.75,7.5);
\draw[thick, dotted] (6,5)--(6,5.5);
\draw (1,0) rectangle (8,10);
\etikz &=& \btikz[line width=.5,scale=.5,baseline=70]

\draw (6,8.5)--(6.5,7.75)--(7,8.5);
\draw (6,1.5)--(6.5,2.25)--(7,1.5);

\draw[white, line width =7.5] (6.5,5)--(6.5,7.5);
\draw (6.5,2.25)--(6.5,7.75);


\begin{scope}[ >=stealth,decoration={markings, mark=at position .6 with {\arrow{>}}}]
\draw[white, line width=7.5] (6.5,5) circle [x radius=1, y radius=.5, start angle = 0, end angle = -180];
\draw[postaction={decorate},blue] (6.5,5) circle [x radius=1, y radius=.5];
\end{scope}

\draw (3,0)--(3,10);
\draw[blue] (3,5)--(4,5);
\draw (3.5,5) node[below] {$\psi$};
\draw (4,0)--(4,10);

\draw (6,4.5) node[below] {$\psi$};

\foreach \y in {6,7}
{
\draw[white,line width=7.5] (\y,0)--(\y,1.5);
\draw[white,line width=7.5] (\y,8.5)--(\y,10);
\draw (\y,0)--(\y,1.5);
\draw (\y,8.5)--(\y,10);
}


\draw[white,line width=7.5] (5,0)--(5,10);
\draw (5,0)--(5,10);
\foreach \x in {1,2,5,8}
{\draw (\x,0)--(\x,10);}

\draw[white, line width =7.5] (6.5,5.25)--(6.5,7.5);
\draw (6.5,5.25)--(6.5,7.75);

\foreach \x in {1,2,3,4}
{\draw (\x,.375) node [left] {\small $\sigma$};}
\foreach \y in {5,6,7,8}
{\draw (\y,.375) node [left] {\small $X$};}
\draw (6.5,2.5) node[right] {$c$};
\stdlabels{\psi}{1}{\psi}{1}
\draw (1,0) rectangle (8,10);

\etikz \\ &= \left ( \frac{S_{\psi 1, c}}{S_{11,c}} \right) & \btikz[line width=.5,scale=.5,baseline=70]

\draw (6,8.5)--(6.5,7.75)--(7,8.5);
\draw (6,1.5)--(6.5,2.25)--(7,1.5);

\draw[white, line width =7.5] (6.5,5)--(6.5,7.5);
\draw (6.5,2.25)--(6.5,7.75);


\draw (3,0)--(3,10);
\draw[blue] (3,5)--(4,5);
\draw (3.5,5) node[below] {$\psi$};
\draw (4,0)--(4,10);

\foreach \y in {6,7}
{
\draw[white,line width=7.5] (\y,0)--(\y,1.5);
\draw[white,line width=7.5] (\y,8.5)--(\y,10);
\draw (\y,0)--(\y,1.5);
\draw (\y,8.5)--(\y,10);
}


\draw[white,line width=7.5] (5,0)--(5,10);
\draw (5,0)--(5,10);
\foreach \x in {1,2,5,8}
{\draw (\x,0)--(\x,10);}

\draw[white, line width =7.5] (6.5,5.25)--(6.5,7.5);
\draw (6.5,5.25)--(6.5,7.75);

\foreach \x in {1,2,3,4}
{\draw (\x,.375) node [left] {\small $\sigma$};}
\foreach \y in {5,6,7,8}
{\draw (\y,.375) node [left] {\small $X$};}
\draw (6.5,5) node[right] {$c$};
\stdlabels{\psi}{1}{\psi}{1}
\draw (1,0) rectangle (8,10);

\etikz
\end{eqnarray}
Although we have draw some horizontal lines here, the edges are understood to have the orientation they inherit from reading the diagram from the bottom up. 

What remains can be resolved using a combination of $F$-moves and bubble-popping relations in a similar manner to the final steps of the calculation for $\langle 1 | T | 1 \rangle$. The contributions from each step in the calculation result in the formula

\begin{eqnarray} \frac{1}{d_{X}}\sum_c \sqrt{d_c} (R^{XX}_c)^2 \left |[F^{XXX}_{X}]_{c,11} \right |^2\left( \sqrt{\frac{d_{X}d_{X}}{d_c}} \right)^2 \Big (  [F^{\sigma 1\sigma 1\sigma 1}_{\sigma 1}]_{11,11}\left( d_{\sigma 1}^2 \sqrt{d_c} d_{X}\right) \\
 + \Big(\frac{S_{\psi 1, c}}{S_{11,c}}\Big ) [F^{\sigma 1\sigma 1\sigma 1}_{\sigma 1}]_{\psi 1,11} [F^{\sigma 1\sigma 1\sigma 1}_{\sigma 1}]_{\psi 1,\psi 1}\left (\sqrt{ \frac{d_{\sigma 1}d_{\sigma 1}}{d_{\psi 1}}}\right )^2 \left(d_{\sigma 1}\sqrt{d_c} d_X \right) \\ = d_{\sigma 1}d_{X_1}^2 \sum_c  \left ( R^{X_1 X_1}_c\right )^2 \left | \left [ F^{X_1X_1X_1}_{X_1}\right ]_{c,11} \right |^2 \left ( 1 - \left ( \frac{S_{\psi 1,c}}{S_{11,c}} \right) \right)
 \end{eqnarray} 

These formulas result in the ratio $\frac{\langle \psi | T | \psi \rangle }{\langle 1 | T | 1 \rangle} = e^{\pi i /4}$.

We stress that the sequence of diagrammatic moves made here is constrained by what is known about the $G$-crossed data for the extension of the bilayer Ising category, and is only one sequence of many possible moves that calculates the same quantity. It is interesting to note that the calculation is independent of $U$- and $\eta$- symbols. 

\end{document}